\documentclass{article}
\usepackage{amsthm,amsmath,amssymb}
\usepackage[margin=1in]{geometry} 
\usepackage{caption}
\usepackage{subcaption}
\usepackage{graphicx}
\usepackage{epstopdf}
\usepackage{color}
\usepackage{framed}
\usepackage{lscape}
\usepackage{rotating}
\usepackage{algorithm}
\usepackage[noend]{algpseudocode}
\usepackage{grffile}
\usepackage{multirow}
\usepackage{xcolor}
\usepackage{comment}
\usepackage{xr}
\usepackage[colorlinks=true,linkcolor=black,citecolor=black,urlcolor=black]{hyperref}

\newtheorem{thm}{Theorem}[section]

\newtheorem{prop}[thm]{Proposition}

\theoremstyle{definition}
\newtheorem{rmk}{Remark}

\newcommand{\EE}{\mathbb{E}}

\newcommand{\tr}[1]{\text{tr}(#1)}

\newcommand{\snr}{\mathrm{SNR}}

\newcommand{\R}{\mathbb{R}}
\newcommand{\A}{\mathcal{A}}
\newcommand{\AMSE}{\text{AMSE}}

\newcommand{\la}{\langle}
\newcommand{\argmin}{\operatorname*{argmin}}
\newcommand{\ra}{\rangle}
\newcommand{\F}{\mathcal{F}}

\newcommand{\naive}{\mathrm{naive}}

\newcommand{\dd}{\mathrm{dd}}

\newcommand{\U}{\mathbf{U}}
\newcommand{\V}{\mathbf{V}}

\renewcommand{\SS}{\mathcal{S}}

\newcommand{\KN}{\mathrm{KN}}

\newcommand{\op}{\mathrm{op}}

\renewcommand{\t}{\mathbf{t}}

\renewcommand{\L}{\mathcal{L}}

\newcommand{\Prob}{\mathbb{P}}

\newcommand{\diag}{\mathrm{diag}}
\newcommand{\what}{\widehat}
\newcommand{\wtilde}{\widetilde}

\newcommand{\loc}{\mathrm{loc}}
\newcommand{\shr}{\mathrm{shr}}

\newcommand{\nuc}{\mathrm{nuc}}

\newcommand{\X}{\mathbf{X}}
\newcommand{\G}{\mathbf{G}}
\renewcommand{\P}{\mathbf{P}}
\newcommand{\Q}{\mathbf{Q}}
\newcommand{\Y}{\mathbf{Y}}
\newcommand{\N}{\mathbf{N}}

\newcommand{\q}{\mathbf{q}}

\newcommand{\vphi}{\varphi}

\renewcommand{\a}{\mathbf{a}}
\renewcommand{\b}{\mathbf{b}}

\renewcommand{\c}{\mathbf{c}}
\renewcommand{\d}{\mathbf{d}}

\renewcommand{\u}{\mathbf{u}}
\renewcommand{\v}{\mathbf{v}}

\newcommand{\x}{\mathbf{x}}
\newcommand{\y}{\mathbf{y}}

\newcommand{\w}{\mathbf{w}}

\newcommand{\I}{\mathbf{I}}
\renewcommand{\A}{\mathbf{A}}
\newcommand{\C}{\mathbf{C}}
\newcommand{\D}{\mathbf{D}}
\newcommand{\E}{\mathbf{E}}

\newcommand{\B}{\mathbf{B}}

\newcommand{\TT}{\mathcal{T}}
\newcommand{\T}{\mathbf{T}}
\renewcommand{\S}{\mathbf{S}}

\newcommand{\W}{\mathbf{W}}
\newcommand{\Z}{\mathbf{Z}}

\newcommand{\Fr}{\mathrm{F}}

\begin{document}
\title{Matrix denoising for weighted loss functions and \\
        heterogeneous signals}
\author{William Leeb\thanks{School of Mathematics, University of Minnesota, Twin Cities. Minneapolis, MN.}}
\date{}
\maketitle

\begin{abstract}
We consider the problem of estimating a low-rank matrix from a noisy observed matrix. Previous work has shown that the optimal method depends crucially on the choice of loss function. In this paper, we use a family of weighted loss functions, which arise naturally for problems such as submatrix denoising, denoising with heteroscedastic noise, and denoising with missing data. However, weighted loss functions are challenging to analyze because they are not orthogonally-invariant. We derive optimal spectral denoisers for these weighted loss functions. By combining different weights, we then use these optimal denoisers to construct a new denoiser that exploits heterogeneity in the signal matrix to boost estimation with unweighted loss.
\end{abstract}

%%%\tableofcontents \newpage

%
%

\section{Introduction}
\label{sec-intro}

This paper is concerned with estimating a low-rank signal matrix $\X$ from an observed matrix $\Y = \X + \G$, where $\G$ is a full-rank matrix of noise. We consider two distinct aspects of the matrix denoising problem. First, we study methods designed for a broader family of loss functions, known as \emph{weighted} loss functions, than considered in earlier works. Second, we design a new denoiser for \emph{unweighted} loss that improves upon previous work by exploiting \emph{heterogeneity} in the target matrix's singular vectors. Like many works on matrix denoising, our methods are designed for an asymptotic regime where the number of rows and columns of $\X$ grow infinitely large, and where the energy in the noise swamps the energy in the signal. This setting is often referred to as the \emph{spiked model} \cite{bai2012sample, bai2008central, baik2006eigenvalues, paul2007asymptotics, johnstone2001distribution, benaych2012singular}.

The methods introduced in this paper extend singular value shrinkage \cite{shabalin2013reconstruction, gavish2017optimal, gavish2014optimal, nadakuditi2014optshrink, donoho2014minimax, leeb2020operator}, which modifies $\Y$'s singular values to mitigate the effects of noise. Our method of \emph{spectral denoising} agrees with singular value shrinkage with unweighted loss, but performs better with weighted loss. While weighted loss functions arise in a number of applications which we describe, they are  challenging as they are not orthogonally-invariant. To derive optimal spectral denoisers for weighted loss, we extend the asymptotic theory of the spiked model, building on work from \cite{leeb2019optimal}.

Our new method of \emph{localized denoising} is designed for unweighted loss. Unlike singular value shrinkage, however, localized denoising exploits heterogeneity in $\X$'s singular vectors;  when certain blocks of coordinates of $\X$ are known to contain more of the signal's energy than others, localized denoising outperforms shrinkage. At the same time, localized denoising's asymptotic performance is never worse than shrinkage's, and so localized denoising inherits shrinkage's  well-known optimality properties.

\subsection{Main ideas}

\begin{table}
\centering
\begin{tabular}{| c | c  | c | }  
\hline  
\# &  Description & Reference  \\
\hline  
\ref{alg:denoise} & Optimal spectral denoising for weighted loss & Section \ref{sec-minfro} \\
\ref{alg-local} & Localized denoising for unweighted loss & Section \ref{sec-localized} \\
\ref{alg-submatrix} & Submatrix denoising & Section \ref{sec-splitting} \\
\ref{alg-hetero} & Matrix denoising with doubly-heteroscedastic noise 
    & Section \ref{sec-whitening} \\
\ref{alg-missing} & Matrix denoising with missing data & Section \ref{sec-missing} \\
\hline
\end{tabular}
\caption{Algorithms introduced in this paper.}
\label{table:algorithms}
\end{table}

In the high-noise, high-dimensional spiked model, the energy of the noise $\G$ is unbounded as $p,n \to \infty$, while the energy of $\X$ is fixed. Consistent estimation of $\X$ from $\Y$ is therefore not possible, so the ``best'' denoiser depends on the choice of loss function. The weighted loss functions we use arise in a variety of applications, described in Section \ref{sec-motivation}; the new method of spectral denoising is adapted to each of these. Table \ref{table:algorithms} lists these algorithms and their locations in the paper.

%%%Because weighted losses are not orthogonally-invariant, deriving the optimal spectral denoiser necessitates the derivation of new asymptotic results on the high-dimensional spiked model.

The optimal spectral denoiser for weighted loss solves a least-squares problem parameterized by weighted inner products between the singular vectors of $\X$ and $\Y$. Though formulas for unweighted inner products  are well-known \cite{paul2007asymptotics, benaych2012singular}, the results we need require a new analysis extending our earlier work in \cite{leeb2019optimal}. While we leave the details to Theorem \ref{thm-asymptotics}, the key idea is that a singular vector $\hat \u_j$ of $\Y$ may be written as a combination of its projection onto the corresponding singular vector $\u_j$ of $\X$ and a residual unit vector $\tilde \u_j$, $\hat \u_j = c_j \u_j + s_j \tilde \u_j$. Here, $c_j$ and $s_j$ are known from the classical theory of the spiked model \cite{paul2007asymptotics}. Because the noise $\G$ is orthogonally invariant, the $\tilde \u_j$ are uniformly random in the subspace orthogonal to $\X$'s singular vectors. Consequently, inner products of the form $\tilde \u_j^T \A \tilde \u_k$ have predictable behavior when the dimension is large \cite{hanson1971bound, wright1973bound, rudelson2013hanson}.

%%%That is, even though the vectors $\tilde \u_j$ are not known, their inner products behave predictably when $p$ is large.

%%%Though the analysis requires orthogonally-invariant noise, we conjecture that the results are more general; indeed, all that is needed is for the quadratic forms $\tilde \u_j^T \A \tilde \u_k$ to concentrate around their expected values. In Section \ref{sec-nongaussian}, we demonstrate in numerical simulations that for certain non-Gaussian distributions with sufficiently many finite moments, the finite sample statistics of the weighted inner products match those for Gaussian noise. By contrast, the results break down for noise with heavy tails.

\subsection{Illustrative example}
\label{sec-intro-localized}
\begin{figure}
\centering
\includegraphics[scale=0.4]{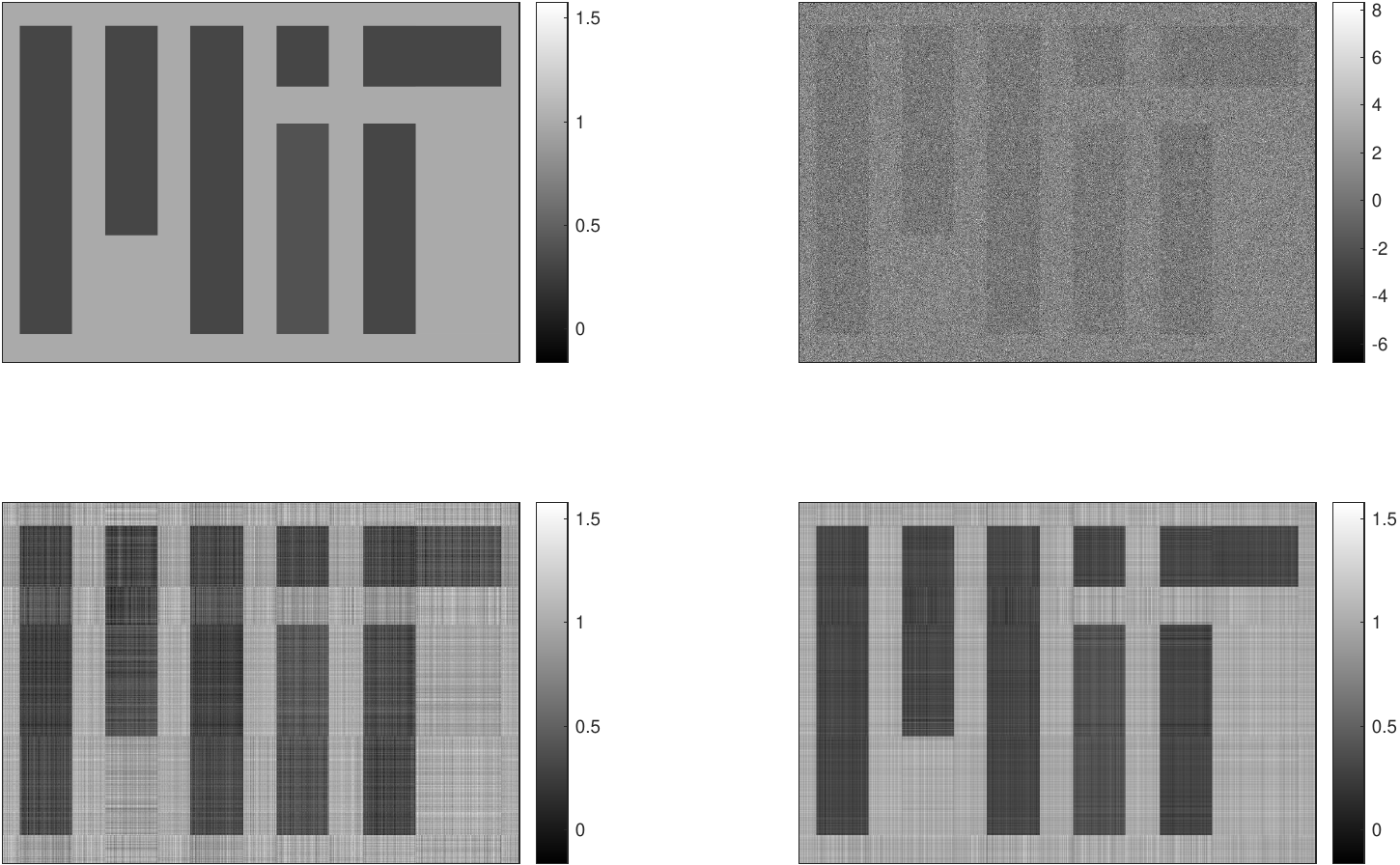}
\caption{
Denoising the MIT logo; see Section \ref{sec-intro-localized} for details. Upper left: the rank $5$ signal matrix. Upper right: the observed noisy matrix. Lower left: the matrix denoised by optimal singular value shrinkage \cite{gavish2017optimal, shabalin2013reconstruction}. Lower right: the matrix denoised by localized denoising (Algorithm \ref{alg-local}). The relative error of singular value shrinkage is approximately $1.25 \times 10^{-1}$, whereas the relative error of localized denoising is approximately $7.41 \times 10^{-2}$.
}
\label{fig-mit}
\end{figure}

The method of \emph{localized denoising}, introduced in Section \ref{sec-localized}, uses the optimal spectral denoiser for \emph{weighted} loss to construct a matrix denoiser for \emph{unweighted} loss. The matrix is broken into submatrices, each of which is denoised by applying the optimal spectral denoiser with weights projecting onto that submatrix's coordinates. In Figure \ref{fig-mit}, we illustrate the performance of localized denoising on the MIT logo, which is a $1574$-by-$2800$ matrix with rank $5$. The logo is corrupted by iid Gaussian noise with standard deviation $\sigma = t_5 / (1.5\gamma^{1/4})$, $t_5$ being the smallest singular value of the clean image. We apply optimal singular value shrinkage and localized denoising, the latter by breaking the rows into $15$ equispaced segments and the columns into $30$ equispaced segments. The relative error $\|\what \X^{\loc} - \X\|_{\Fr} / \|\X\|_{\Fr}$ of localized denoising is approximately $7.41 \times 10^{-2}$; the relative error of singular value shrinkage is approximately $1.25 \times 10^{-1}$, which is significantly larger. The improvement from localized denoising is due to the signal matrix's heterogeneity along the rows and columns. However, the row and column subdivisions are not chosen to extract any specific structure in the image, and localized denoising does not appear to be very sensitive to the choice of subdivisions; similar results may be obtained with other subdivisions as well.

%%% The locally denoised image is also visually much closer to the true logo.

%%%; the dark and light regions are more homogeneous.

%%%; the choice of $15$ and $30$ is arbitrary, and other choices give approximately the same performance.

\subsection{Outline of the paper}

Section \ref{sec-prelim} contains the problem statement and key definitions. Section \ref{sec-math} presents the new asymptotic results. Section \ref{sec-minfro} derives the optimal spectral denoiser for weighted loss. Section \ref{sec-localized} introduces localized denoising. Section \ref{sec-motivation} describes three applications of weighted loss functions. Section \ref{sec-numerics} reports on numerical experiments. Section \ref{sec-conclusion} concludes by discussing potential applications.

\section{Preliminaries}
\label{sec-prelim}

\subsection{The observation model}
\label{sec-model}

We observe a $p$-by-$n$ data matrix $\Y = \X+\G$, consisting of a low-rank signal matrix $\X$ and a full-rank isotropic Gaussian noise matrix $\G$. We write $\X$ as
\begin{math}
\X = \sum_{k=1}^{r} t_k \u_k \v_k^T,
\end{math}
where the $\u_k$ and $\v_k$ are orthonormal vectors in $\R^p$ and $\R^n$, respectively, and $t_1 > \dots > t_r > 0$. The entries of the noise matrix $\G$ are iid $N(0,1/n)$. We write $\Y$ as
\begin{math}
\Y = \sum_{k=1}^{\min(n,p)} \lambda_k \hat \u_k \hat \v_k^T,
\end{math}
where the $\hat \u_k$ and $\hat \v_k$ are orthonormal vectors in $\R^p$ and $\R^n$, respectively, and $\lambda_1  \ge \dots \ge \lambda_{\min(n,p)} \ge 0$.

We let $\Omega = \Omega_p$ be one of a sequence of matrices with $p$ columns, and $\Pi = \Pi_n$ be one of a sequence of matrices with $n$ columns. In Section \ref{sec-loss}, these matrices will be used to define the loss function for estimating $\X$. In order to have a well-defined asymptotic theory when $p, n \to \infty$, we assume that certain quantities defined in terms of $\Omega$, $\Pi$, and the singular vectors of $\X$ have definite, finite limits. We define:
\begin{align}
\label{def-mu}
\mu = \lim_{p \to \infty} \frac{1}{p} \tr{\Omega^T \Omega},
\quad
\nu = \lim_{n \to \infty} \frac{1}{n} \tr{\Pi^T \Pi}.
\end{align}
For $1 \le j, k \le r$, we assume that the weighted inner products between the population singular vectors converge almost surely in the large $p$, large $n$ limits:
\begin{align}
\label{eq-ejk}
e_{jk} = \lim_{p \to \infty} \la \Omega \u_j, \Omega \u_k\ra,
\quad
\tilde e_{jk} = \lim_{n \to \infty} \la\Pi \v_j, \Pi \v_k\ra.
\end{align}
For $1 \le k \le r$, we will let $\alpha_k = e_{kk}$ and $\beta_k = \tilde e_{kk}$.

We assume that these limits exist and are finite and positive.  We also assume that the operator norms of the matrices $\Omega = \Omega_p$ and $\Pi = \Pi_n$  remain bounded as $p,n \to \infty$; that is, $\|\Omega_p \|_{\op}, \|\Pi_n\|_{\op} \le C < \infty$ for all $p$ and $n$, where $C$ does not depend on $p$ or $n$. These conditions are the only assumptions we make on the matrices $\Omega$ and $\Pi$.

We will parametrize the problem size by the number of columns $n$, and let the number of rows $p = p_n$ grow with $n$. Specifically, we will assume that the limit
\begin{align}
\gamma = \lim_{n \to \infty} \frac{p_n}{n}
\end{align}
is well-defined and finite. In all statements where $n \to \infty$, it will be implicitly assumed as well that $p \to \infty$ and $p/n \to \gamma$. We assume that the number of population components $r$ and the singular values $t_1,\dots,t_r$ stay fixed with $p$ and $n$.

\begin{rmk}
\label{rmk-notation}
All quantities that depend on $p$ and $n$, such as $\u_k$ and $\v_k$, are actually elements of a \emph{sequence} indexed by $p$ and/or $n$. However, for notational simplicity, we drop the explicit dependence on $p$ and $n$ unless it is needed for clarity.
\end{rmk}

\begin{rmk}
There are counterexamples to the existence of the limits in \eqref{def-mu} and \eqref{eq-ejk}. For example, one may take $\u_1$ to be the first standard unit vector $(1,0,\dots,0)^T$ when $p$ is even, and the constant vector $(1,\dots,1)^T / \sqrt{p}$ when $p$ is odd; and take $\Omega_p = \diag(0,1,\dots,1)$. The limit defining $\alpha_1$ will not exist in this case, as odd terms in the sequence $\|\Omega_p \u_1\|^2$ converge to $1$, and even terms converge to $0$. By contrast, the examples in Section \ref{sec-numerics} satisfy the asymptotic conditions.
\end{rmk}

\begin{rmk}
The values $\mu$, $\nu$, $e_{jk}$, and $\tilde e_{jk}$ from equations \eqref{def-mu} and \eqref{eq-ejk} are used to characterize the weighted inner products between singular vectors of $\X$ and $\Y$; see Theorem \ref{thm-asymptotics}. These weighted inner products are needed to evaluate optimal spectral denoisers for weighted loss, as described in Section \ref{sec-minfro}.
\end{rmk}

\subsection{Heterogeneity, genericity, and weighted orthogonality}
\label{sec-het}

One of the aspects of the theory of matrix denoising we will explore is the role of the signal matrix $\X$'s singular vectors, $\u_1,\dots,\u_r$ and $\v_1,\dots,\v_r$. To that end, we introduce two definitions we will be using throughout the paper. We say that a unit vector $\x \in \R^m$ is \emph{generic} with respect to an $m$-by-$m$ positive-semidefinite matrix $\A_m \in \R^{m \times m}$ if
\begin{math}
\x^T \A \x \sim \frac{1}{m} \tr{\A},
\end{math}
where ``$\sim$'' indicates that the difference between the two sides vanishes almost surely as $m \to \infty$ (to be precise, $\x$ and $\A$ are elements of a sequence of vectors and matrices, respectively, indexed by $m$; but following the convention described in Remark \ref{rmk-notation} we will drop the explicit dependence on $m$).

By contrast, we say that $\x$ is \emph{heterogeneous} if it is not generic. This means that the energy of the vector $\x$ is not uniformly distributed across its coordinates in the eigenbasis of $\A$. Indeed, if $\A = \sum_{k=1}^{m} h_k \w_k \w_k^T$ is the eigendecomposition of $\A$, then 
\begin{align}
\x^T \A \x = \sum_{k=1}^{m} h_k \la \x,\w_k \ra^2.
\end{align}
If the energy of $\x$ were equally spread out across the $\w_k$, then $\la \x , \w_k \ra \sim 1/\sqrt{m}$, and so $\x^T \A \x \sim \tr{\A} / m$.

Given a collection of vectors $\x_1,\dots,\x_k \in \R^{m}$, we will say that they satisfy the \emph{weighted orthogonality condition} (or are \emph{weighted orthogonal}) with respect to a positive-semidefinite matrix $\A$ if 
\begin{align}
\x_i^T \A \x_j \sim 0
\end{align}
whenever $i \ne j$. In other words, the $\x_j$ are asymptotically orthogonal with respect to the weighted inner product defined by $\A$. 

\begin{rmk}
From the Hanson-Wright inequality \cite{hanson1971bound, wright1973bound, rudelson2013hanson}, random unit vectors $\x$ from suitably regular distributions are generic, with respect to any $\A$ with bounded operator norm. Furthermore, the weighted orthogonality condition will also hold for independent random unit vectors $\x_1,\dots,\x_k$ from a suitable distribution (see \cite{benaych2011fluctuations}).
\end{rmk}

%%%\subsection{Spectral denoisers}

\subsection{Spectral denoisers and weighted loss functions}

\label{sec-loss}

For the top $r$ empirical singular vectors $\hat \u_1, \dots, \hat \u_r$ and $\hat \v_1,\dots, \hat \v_r$ of $\Y$, define the matrices $\what \U = [\hat \u_1,\dots, \hat \u_r]$ and $\what \V = [\hat \v_1,\dots, \hat \v_r]$. Consider the class of estimators  defined by
\begin{align}
\SS = \left\{ \what \U \what \B \what \V^T : \what \B \in \R^{r \times r} \right\}.
\end{align}
Each matrix in $\SS$ has the same singular subspaces as the observed matrix $\Y$, though not necessarily the same singular vectors. We call $\SS$ the family of \emph{spectral denoisers}.

%

%%%\subsection{The choice of loss function}
%%%\label{sec-loss}

%%%In the high-dimensional regime, where $p$ grows proportionally to $n$, it is known that the method for estimating $\X$ depends crucially on the choice of loss function \cite{gavish2017optimal}. 
%
%%%Up until now, the only loss functions that have been considered are orthogonally-invariant. In the case of an isotropic Gaussian noise matrix (as we consider here), the paper \cite{gavish2017optimal} derives optimal shrinkers for loss-functions that are orthogonally-invariant and block-decomposable, meaning that if both $\what \X$ and $\X$ have the same block-diagonal structure, then the loss can be written as functions of the losses between the individual blocks.
%
%%%In this work, we will consider a different class of loss functions than have been previously considered.

We consider estimating the low-rank signal matrix $\X$ with respect to the \emph{weighted Frobenius loss} defined by
\begin{align}
\L_n(\what \X, \X) = \|\Omega (\what \X -  \X ) \Pi^T \|_{\Fr}^2,
\end{align}
where $\|\cdot\|_{\Fr}$ denotes the matrix Frobenius norm, and $\Omega$ and $\Pi$ are matrices satisfying the conditions in Section \ref{sec-model}. This type of loss function is used when the user pays different prices for errors in different rows and columns.

We now define the precise estimation problem we will consider. For any deterministic $r$-by-$r$ matrix $\what \B$, we define the asymptotic error
\begin{align}
\label{eq-limitloss}
\L(\what \U \what \B \what \V^T, \X) = \lim_{n \to \infty} \L_n(\what \U \what \B \what \V^T, \X).
\end{align}
Our goal is then to find the matrix $\what \B$ to minimize this loss, and show how $\what \B$ may be consistently estimated from the observed matrix $\Y$. That is, we define
\begin{align}
\label{eq-optB}
\what \B 
= \argmin_{\what \B' \in \R^{r \times r}}
    \L(\what \U \what \B' \what \V^T, \X)
\end{align}
and define $\what \X = \what \U \what \B \what \V^T$.

\begin{rmk}
For any deterministic $\what \B$, the asymptotic loss \eqref{eq-limitloss} exists and is finite almost surely, even though the matrices $\what \U \what \B \what \V^T$ and $\X$ are growing in size. It will be shown in Section \ref{sec-minfro} that since $\what \U \what \B \what \V^T$ and $\X$ each have rank at most $r$, $\|\Omega (\what \X -  \X ) \Pi^T \|_{\Fr}^2$ depends only on $t_1,\dots,t_r$; the $r^2$ entries of $\what \B$; and the weighted inner products between the top $r$ singular vectors of $\Y$ and $\X$. It will follow from Theorem \ref{thm-asymptotics} that these inner products converge almost surely to finite limits, and consequently that the asymptotic loss \eqref{eq-limitloss} is well-defined almost surely.
%%%Analogous observations are used to define the asymptotic loss for singular value shrinkage and eigenvalue shrinkage with unweighted loss \cite{gavish2017optimal, nadakuditi2014optshrink, donoho2018optimal}.
\end{rmk}

\section{Asymptotic theory for the spiked model}
\label{sec-math}

In this section, we derive the limits of inner products between the weighted population and empirical vectors. We define the cosines between the unweighted empirical and population vectors:
\begin{align}
\label{cosines-unweighted}
c_{jk} = \lim_{p \to \infty} \la \hat \u_j, \u_k \ra,
    \quad
\tilde c_{jk} = \lim_{n \to \infty} \la \hat \v_j, \v_k \ra.
\end{align}
Next we define \emph{weighted} inner products between the population and empirical vectors:
\begin{align}
\label{eq-cjk}
c_{jk}^\omega = \lim_{p \to \infty} \la \Omega \hat \u_j, \Omega \u_k \ra,
    \quad
\tilde c_{jk}^\omega = \lim_{n \to \infty} \la \Pi \hat \v_j, \Pi \v_k\ra.
    \quad
\end{align}
These are inner products with weight matrices $\Omega^T \Omega$ and $\Pi^T \Pi$, respectively. We also define the weighted inner products between the empirical singular vectors:
\begin{align}
\label{eq-djk}
d_{jk} = \lim_{p \to \infty} \la \Omega \hat \u_j, \Omega \hat \u_k\ra,
    \quad
\tilde d_{jk} = \lim_{n \to \infty} \la \Pi \hat \v_j, \Pi \hat \v_k\ra.
\end{align}
We let $c_k^\omega = c_{kk}^\omega$ and $\tilde c_k^\omega = \tilde c_{kk}^\omega$, and similarly for the other terms.

\begin{rmk}
From Theorem \ref{thm-asymptotics} below, the limits \eqref{eq-cjk} and \eqref{eq-djk} exist almost surely and are finite so long as the assumptions on $\Omega$ and $\Pi$ from Section \ref{sec-model} hold.
\end{rmk}

The first result provides formulas for $c_{jk}$ and $\tilde c_{jk}$, and relates the singular values of $\X$ to those of $\Y$. It is well-known in the literature (see e.g.\ \cite{paul2007asymptotics, benaych2012singular}).

\begin{prop}
\label{prop-vanilla}

For $1 \le k \le r$, the $k^{th}$ singular value of $\Y$ converges almost surely as $p,n \to \infty$ to $\lambda_k$, defined by:
\begin{align}
\lambda_k^2= 
\begin{cases}
(t_k^2 + 1)\left( 1 + \frac{\gamma}{t_k^2}\right),
        & \text{ if } t_k  > \gamma^{1/4}, \\
(1 + \sqrt{\gamma})^2, & \text{ if }  t_k \le \gamma^{1/4},
\end{cases}
\end{align}

For $1 \le j, k \le r$, the limits \eqref{cosines-unweighted} defining $c_{jk}$ and $\tilde c_{jk}$ almost surely exist and are given by the following expressions:
\begin{align}
c_{jk}^2 = 
\begin{cases}
\frac{1 - \gamma / t_k^4}{1 + \gamma / t_k^2},
        & \text{ if } j = k \text{ and } t_k  > \gamma^{1/4}, \\
0, & \text{ if } j \ne k \text{ or } t_k \le \gamma^{1/4},
\end{cases}
\end{align}
and
\begin{align}
\tilde c_{jk}^2 = 
\begin{cases}
\frac{1 - \gamma / t_k^4}{1 + 1 / t_k^2 },
        & \text{ if } j = k \text{ and } t_k > \gamma^{1/4}, \\
0, & \text{ if } j \ne k \text{ or } t_k \le \gamma^{1/4}.
\end{cases}
\end{align}

\end{prop}

\begin{rmk}
\label{rmk-positive}
While the signs of $c_{k}$ and $\tilde c_{k}$ are arbitrary, their product satisfies $c_{k} \tilde c_{k} \ge 0$. We may therefore assume that $c_{k} \ge 0$ and $\tilde c_{k} \ge 0$
(see, e.g., \cite{nadakuditi2014optshrink}).
\end{rmk}

%%%The next result describes the limits of the weighted inner products.

\begin{thm}
\label{thm-asymptotics}
Suppose $1 \le j, k \le r$. Then the limits \eqref{eq-cjk} and \eqref{eq-djk} almost surely exist and are equal to the following expressions:
\begin{align}
\label{fmla-comega}
c_{jk}^\omega = 
\begin{cases}
e_{jk} c_j, & \text{ if }  t_j > \gamma^{1/4}, \\
0, & \text{ if }  t_j \le \gamma^{1/4},
\end{cases}
\end{align}
\begin{align}
\tilde c_{jk}^\omega = 
\begin{cases}
\tilde e_{jk} \tilde c_j,  & \text{ if }  t_j > \gamma^{1/4}, \\
0, & \text{ if }  t_j \le \gamma^{1/4},
\end{cases}
\end{align}
\begin{align}
\label{fmla-d}
d_{jk} = 
\begin{cases}
c_k^2 \alpha_k + s_k^2 \mu,
        & \text{ if }  j = k \text{ and } t_k > \gamma^{1/4}, \\
e_{jk} c_j c_k ,
        & \text{ if }  j \ne k \text{ and } \min\{t_j,t_k\} > \gamma^{1/4}, \\
0, & \text{ if }  j \ne k \text{ and } \min\{t_j,t_k\} \le \gamma^{1/4},
\end{cases}
\end{align}
\begin{align}
\tilde d_{jk} = 
\begin{cases}
\tilde c_k^2 \beta_k + \tilde s_k^2 \nu,
        & \text{ if }  j = k \text{ and } t_k > \gamma^{1/4}, \\
\tilde e_{jk} \tilde c_j \tilde c_k ,
        & \text{ if }  j \ne k \text{ and } \min\{t_j,t_k\} > \gamma^{1/4}, \\
0, & \text{ if }  j \ne k \text{ and } \min\{t_j,t_k\} \le \gamma^{1/4}.
\end{cases}
\end{align}

\end{thm}

The proof of Theorem \ref{thm-asymptotics} may be found in Section \ref{appendix-asymptotics}.

\begin{rmk}
While the signs of inner products between singular vectors are arbitrary, Theorem \ref{thm-asymptotics} states that once the signs of $e_{jk}$ and $\tilde e_{jk}$ are fixed, the signs of $c_{jk}^\omega$, $\tilde c_{jk}^\omega$, $d_{jk}$ and $\tilde d_{jk}$ are determined.
\end{rmk}

%

%%%\subsection{Minimizing Frobenius loss}

\section{Optimal spectral denoising}
\label{sec-minfro}

In this section, we derive the asymptotically optimal spectral denoiser with respect to the weighted loss \eqref{eq-limitloss}, and show how to consistently estimate it from $\Y$. We define the $r$-by-$r$ weighted inner product matrices $\D = (d_{kl})$, $\wtilde \D = (\tilde d_{jk})$, $\E = (e_{jk})$, $\wtilde \E = (\tilde e_{jk})$, $\C = (c_{jk}^\omega)$, and $\wtilde \C = (\tilde c_{jk}^\omega)$, and the vector $\t = (t_1,\dots,t_r)^T$ of population singular values.
\begin{thm}
\label{thm-ls}
The optimal choice of $\what \B$ is given by
\begin{math}
\what \B = \D^+ \C \diag(\t) \wtilde \C^T \wtilde \D^+,
\end{math}
with weighted AMSE almost surely equal to
\begin{align}
\lim_{n \to \infty} \| \Omega (\what \X - \X) \Pi^T\|_{\Fr}^2
= \la \E \diag(\t) \wtilde \E - 
    \C^T \D^+ \C  \diag(\t) \wtilde \C^T \wtilde \D^+ \wtilde\C,\diag(\t)\ra_{\Fr}.
\end{align}
\end{thm}
The proof of Theorem \ref{thm-ls} may be found in Section \ref{appendix-ls}.

\begin{algorithm}[ht]
\caption{Optimal spectral denoising for weighted loss}
\label{alg:denoise}
\begin{algorithmic}[1]
\item Input: $\Y$; weights $\Omega$ and $\Pi$

\item
rank $r$ SVD of $\Y$:
\begin{tabbing}
\hspace*{0cm}\=\kill
    \> $\lambda_1 \ge \dots \ge \lambda_r > 1+ \sqrt{\gamma}$ \\
    \> $\what \U=[\hat \u_1,\dots,\hat \u_r]$, 
        $\what\V=[\hat \v_1,\dots,\hat \v_r]$
\end{tabbing}

\item
$\mu = \tr{\Omega^T\Omega} / p$,
$\nu = \tr{\Pi^T\Pi} / n$

\item
for $1 \le k \le r$:
\begin{tabbing}
\hspace*{0cm}\=\kill
    \>$t_k = \sqrt{
        \frac{\lambda_k^2 - 1 - \gamma + \sqrt{(\lambda_k^2 - 1 - \gamma)^2 - 4\gamma}}{2}}$
        \\
    \>$c_{k} = \sqrt{ \frac{1 - \gamma / t_k^4}{1 + \gamma / t_k^2} }$,
        $\tilde c_{k} = \sqrt{ \frac{1 - \gamma / t_k^4}{1 + 1 / t_k^2 } }$\\
    \> $s_k = \sqrt{1 - c_k^2}, \tilde s_k = \sqrt{1 - \tilde c_k^2}$\\
    \>$d_{k} = \| \Omega \hat \u_k \|^2$, 
        $\tilde d_{k} = \| \Pi \hat \v_k \|^2$  \\
    \>$\alpha_k = (d_k - s_k^2 \mu) / c_k^2$, 
        $\beta_k = (\tilde d_k - \tilde s_k^2 \nu) / \tilde c_k^2$ \\
    \>$c_{k}^\omega = \alpha_k c_k$, $\tilde c_{k}^\omega = \beta_k \tilde c_k$
\end{tabbing}

\item
for $1 \le j \ne k \le r$:
\begin{tabbing}
\hspace*{0cm}\=\kill
    \>$d_{jk} = \hat \u_j^T \Omega^T \Omega \hat \u_k$, 
        $\tilde d_{jk} = \hat \v_j^T \Pi^T \Pi \hat \v_k$  \\
    \>$e_{jk} = d_{jk}/ (c_j c_k)$, 
         $\tilde e_{jk} = \tilde d_{jk}/ (\tilde c_j \tilde c_k)$,
         $j \ne k$\\
    \>$c_{jk}^\omega = e_{jk} c_j$, $\tilde c_{jk}^\omega = \tilde e_{jk} \tilde c_j$
\end{tabbing}

\item
\begin{tabbing}
\hspace*{0cm}\=\kill
\>$\D = (d_{jk})$, $\wtilde \D = (\tilde d_{jk})$ \\
\>$\E = (e_{jk})$, $\wtilde \E = (\tilde e_{jk})$ \\
\> $\C = (c_{jk}^\omega)$, $\wtilde \C = (\tilde c_{jk}^\omega)$
\end{tabbing}

\item
$\t = (t_1,\dots,t_r)^T$ 

\item
$\what \B = \D^+ \C \diag(\t) \wtilde \C^T \wtilde \D^+$

\item 
$\what \X = \what \U \what \B \what \V^T$

\item
$\AMSE=\la \E \diag(\t) \wtilde \E - 
    \C^T \D^+ \C  \diag(\t) \wtilde \C^T \wtilde \D^+ \wtilde\C,\diag(\t)\ra_{\Fr}$

\end{algorithmic}
\end{algorithm}

The matrices $\D$, $\wtilde \D$, $\E$, $\wtilde \E$, $\C$, and $\wtilde \C$ and the singular values $t_1,\dots,t_r$ may be estimated using Proposition \ref{prop-vanilla} and Theorem \ref{thm-asymptotics}. First, from Proposition \ref{prop-vanilla}, we can estimate $t_k$, $c_k$ and $\tilde c_k$, so long as $t_k > \gamma^{1/4}$, i.e.\ if $\lambda_k > 1 + \sqrt{\gamma}$:
\begin{align}
\label{est-tk}
t_k = \sqrt{
        \frac{\lambda_k^2 - 1 - \gamma + \sqrt{(\lambda_k^2 - 1 - \gamma)^2 - 4\gamma}}{2}
        },
\quad
c_{k} = 
\sqrt{ \frac{1 - \gamma / t_k^4}{1 + \gamma / t_k^2} },
\quad
\tilde c_{k} = 
\sqrt{ \frac{1 - \gamma / t_k^4}{1 + 1 / t_k^2 } }.
\end{align}

\begin{rmk}
From Remark \ref{rmk-positive}, we can take both $c_k$ and $\tilde c_k$ to be positive.
\end{rmk}

The values $d_{jk} = \hat \u_j^T \Omega^T \Omega \hat \u_j$ and $\tilde d_{jk} = \hat \v_j^T \Pi^T \Pi \hat \v_k$ are directly estimable, as they are the weighted inner products between the empirical singular vectors. We then estimate $\alpha_k$ and $\beta_k$, assuming $t_k > \gamma^{1/4}$:
\begin{align}
\label{est-alphak}
\alpha_k = \frac{d_k - s_k^2 \mu}{c_k^2},
\quad
\beta_k = \frac{\tilde d_k - \tilde s_k^2 \nu}{\tilde c_k^2}.
\end{align}

When $j \ne k$, we take $e_{jk} = d_{jk} / (c_j c_k)$ and $\tilde e_{jk} = \tilde d_{jk} / (\tilde c_j \tilde c_k)$ (so long as $t_j$ and $t_k$ both exceed $\gamma^{1/4}$, i.e. $\lambda_j$ and $\lambda_k$ both exceed $1+\sqrt{\gamma}$). Finally, for all $j, k$, we take $c_{jk}^\omega = e_{jk} c_j$ and $\tilde c_{jk} = \tilde e_{jk} \tilde c_j$. The method is summarized in Algorithm \ref{alg:denoise}.

\subsection{Diagonal denoisers}
\label{sec-diagonal}
In this section, we consider a subset of spectral denoisers, in which the matrix $\what \B$ is required to be diagonal. More precisely, we search for a vector $\hat \t = (\hat t_1,\dots,\hat t_r)^T$ of real numbers, so that the estimator
\begin{align}
\what \X^{\mathrm{dd}} = \what \U \diag(\hat \t) \what \V^T
= \sum_{k=1}^{r} \hat t_k \hat \u_k \hat \v_k^T
\end{align}
minimizes the $\AMSE$ $\L(\what \X^{\dd},\X) = \lim_{n \to \infty} \| \Omega (\what \X^{\mathrm{dd}} - \X) \Pi^T\|_{\Fr}^2$.

%%%The $\what \X^{\dd}$ will not outperform the optimal spectral denoiser $\what \X$

\begin{rmk}
Optimal diagonal denoising cannot have better asymptotic performance than optimal spectral denoising, as the diagonal denoiser is a spectral denoiser. However,
Theorem \ref{thm-diag} below shows that under weighted orthogonality, the methods coincide; and the simplicity of $\what \X^{\mathrm{dd}}$ makes it easier to analyze, which will be exploited in the proofs of Theorem \ref{thm-compare2} and Proposition \ref{prop-merge} and the analysis of Section \ref{sec-topt}.
\end{rmk}

%%%We also remark that the weighted orthogonality condition can be tested from the observed matrix $\Y$ using Theorem \ref{thm-asymptotics}, as $e_{jk} = 0$ if and only if $d_{jk} = 0$.

%%%We will show the following result:

\begin{thm}
\label{thm-diag}
Suppose that either $\u_1,\dots,\u_r$ are weighted orthogonal with respect to $\Omega^T\Omega$, or $\v_1,\dots,\v_r$ are weighted orthogonal with respect to $\Pi^T\Pi$. Suppose too that $t_k > \gamma^{1/4}$, $1 \le k \le r$. Then the singular values $\hat t_k$, $1 \le k \le r$, of $\what \X^{\dd}$ are:
\begin{align}
\label{t-hat}
\hat t_k = t_k c_k \tilde c_k
        \cdot \frac{\alpha_k}{c_k^2\alpha_k + s_k^2 \mu}
        \cdot \frac{ \beta_k}{\tilde c_k^2\beta_k + \tilde s_k^2 \nu},
\end{align}
where $t_k$, $c_k$ and $\tilde c_k$ are given by  \eqref{est-tk}, and $\alpha_k$ and $\beta_k$ are given by \eqref{est-alphak}. The weighted AMSE is almost surely equal to
\begin{align}
\lim_{n \to \infty} \| \Omega (\what \X - \X) \Pi^T\|_{\Fr}^2
= \sum_{k=1}^{r} t_k^2 \alpha_k \beta_k \left( 1 
    - c_k^2 \tilde c_k^2 \cdot \frac{\alpha_k}{c_k^2\alpha_k + s_k^2 \mu}
    \cdot \frac{ \beta_k}{\tilde c_k^2\beta_k + \tilde s_k^2 \nu}\right).
\end{align}

If $\u_1,\dots,\u_r$ and $\v_1,\dots,\v_r$ are both weighted orthogonal with respect to $\Omega^T \Omega$ and $\Pi^T \Pi$, respectively, then $\what \X = \what \X^{\dd}$.
%%%That is, the optimal spectral denoiser is equal to the optimal diagonal denoiser.
\end{thm}

The proof of Theorem \ref{thm-diag} is found in Section \ref{appendix-diag}.

%
%    figure from experiments/plot_topt
%
\begin{figure}[h]
\centering
\includegraphics[scale=0.32]{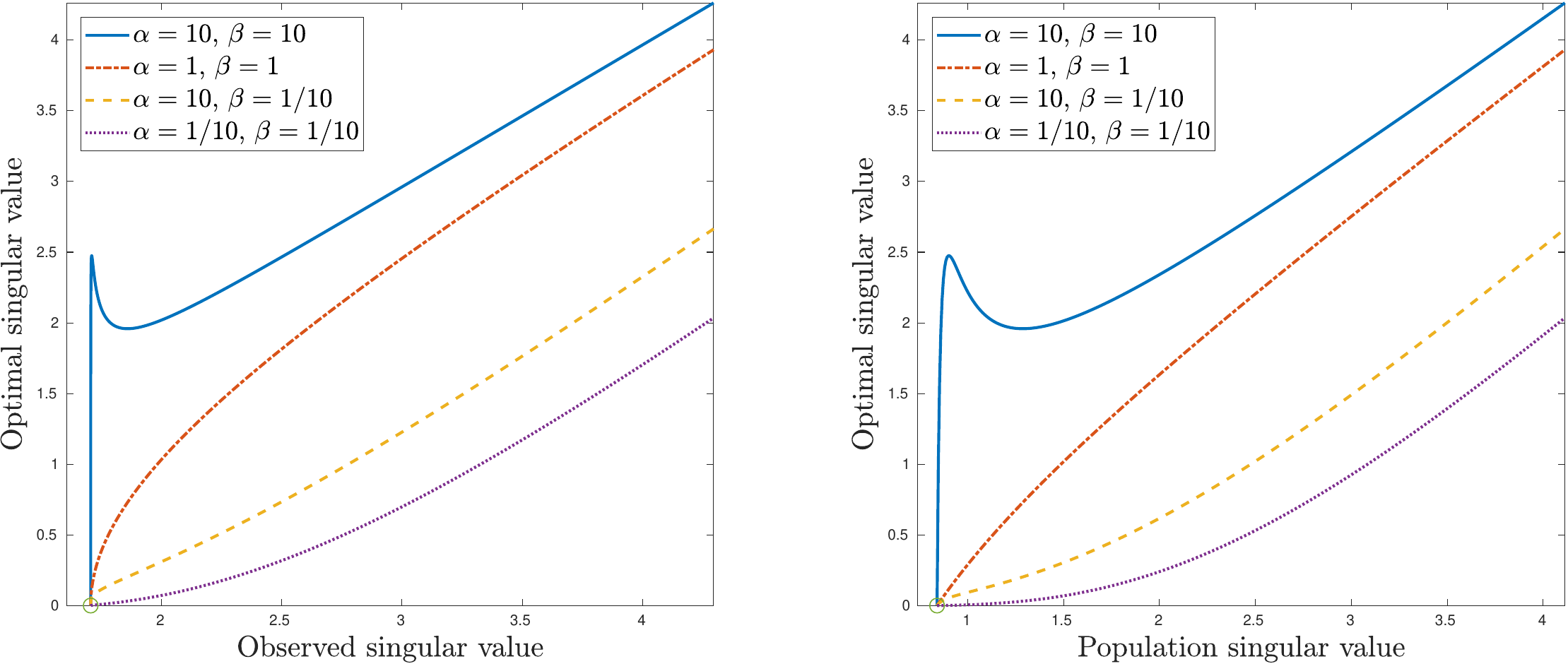}
\caption{
The optimal singular value $\hat t$, plotted as a function of the observed singular value $\lambda$ (left) and the population singular value $t$ (right), for varying values of $\alpha$ and $\beta$ and $\mu = \nu = 1$.
}
\label{fig-shrinkers}
\end{figure}

%

%%%\subsection{Magnitude of $\hat t$}
\subsection{Behavior of the optimal singular values}
\label{sec-topt}

In this section, we assume either that $r=1$; or that $\u_1,\dots,\u_r$ are weighted orthogonal with respect to $\Omega^T\Omega$ and $\v_1,\dots,\v_r$ are weighted orthogonal with respect to $\Pi^T\Pi$. In either case, the optimal spectral denoiser coincides with the optimal diagonal denoiser, and both are given by Theorem \ref{thm-diag}. Though this setting is quite restrictive, it permits us to exploit formula \eqref{t-hat} for the optimal singular values to gain insight into the behavior of the optimal spectral denoiser. Propositions \ref{prop-shrink} and \ref{prop-monotone} describe the behavior of the optimal singular value $\hat t_k$ in this setting. Because each $\hat t_k$ depends only on the information specific to component $k$, we will drop the subscript $k$ from the notation.

\begin{prop}
\label{prop-shrink}
If either $\alpha \le \mu$ or $\beta \le \nu$, then $\hat t \le \lambda$. Conversely, for any fixed value of $t$, there are sufficiently large values of $\alpha$ and $\beta$ for which $\hat t > \lambda$.
\end{prop}

\begin{prop}
\label{prop-monotone}
If $\alpha \le \mu$ or $\beta \le \nu$, then $\hat t$ is an increasing function of $\lambda$.
\end{prop}

The proofs of Propositions \ref{prop-shrink} and \ref{prop-monotone} may be found in Section \ref{appendix-shrink} and Section \ref{appendix-monotone}, respectively.

\begin{rmk}
From \cite{shabalin2013reconstruction, gavish2017optimal, nadakuditi2014optshrink}, the optimal singular value for unweighted Frobenius loss is $\hat t^{\shr} = t c \tilde c$, which is smaller than the observed singular value $\lambda$. Proposition \ref{prop-shrink} shows that with \emph{weighted} loss, such shrinkage only occurs for small $\alpha$ or $\beta$.
\end{rmk}

The conclusion of Proposition \ref{prop-monotone} need not hold if $\alpha > \mu$ and $\beta > \mu$. In Figure \ref{fig-shrinkers} we plot the optimal $\hat t$, both as a function of the observed singular value $\lambda$ and the population singular value $t$, for various values of $\alpha$ and $\beta$ (and $\mu = \nu = 1$). The non-monotonicity is apparent when $\alpha = \beta = 10$.

\section{Localized denoising}
\label{sec-localized}

In this section, we introduce a new procedure called \emph{localized denoising} for estimating $\X$ with \emph{unweighted} Frobenius loss. As we will show, localized denoising is asymptotically never worse than optimal singular value shrinkage \cite{gavish2017optimal, shabalin2013reconstruction}, defined by
\begin{math}
\what \X^{\shr} = \sum_{k=1}^{r} \hat t_k^{\shr} \hat \u_k \hat \v_k^T,
\end{math}
where $\hat t_k^{\shr} = t_k c_k \tilde c_k$. Since singular value shrinkage is optimal for unweighted loss both in the minimax sense and when averaging over a uniform prior on $\u_k$ and $\v_k$ \cite{donoho2014minimax,shabalin2013reconstruction}, localized denoising inherits these same optimality properties. Furthermore, localized denoising can outperform singular value shrinkage when the singular vectors of $\X$ are heterogeneous.

\subsection{Definition of localized denoising}

To define localized denoising, we expand the identity matrices 
\begin{math}
\I_p = \sum_{i=1}^{I} \Omega_i
\end{math}
and
\begin{math}
\I_n = \sum_{j=1}^{J} \Pi_j
\end{math}
into sums of pairwise orthogonal projections $\Omega_i \in \R^{p \times p}$ and $\Pi_j \in \R^{n \times n}$, where $I$ and $J$ are fixed. We require that $\Omega_i = \Omega_i^T = \Omega_i^2$ and $\Omega_{i'} \Omega_i = \mathbf{O}_{p \times p}$ for $i \ne i'$; and similarly for the $\Pi_j$.

%%%Let's describe a particular weighting procedure. Partition the indices $\{1,\dots,p\}$ into $I$ disjoint subsets, $S_1,\dots,S_I$; and partition $\{1,\dots,n\}$ into $J$ disjoint subsets, $T_1,\dots,T_J$. Let $\Omega_i$ be the $p$-by-$p$ diagonal matrix that projects onto the coordinates in $S_i$, and let $\Pi_j$ be the $n$-by-$n$ diagonal matrix that projects onto the coordinates in $T_j$.

%%%More generally, $\Omega_1,\dots,\Omega_I$ may be any projection matrices that sum to the identity, $\sum_{i=1}^{I} \Omega_i = \I_p$; and similarly, $\Pi_1,\dots,\Pi_J$ may be any projection matrices that sum to the identity, $\sum_{j=1}^{J} \Pi_i = \I_n$.

We let $\what \X_{(i,j)}^{\loc}$ denote the optimal spectral denoiser with respect to the weight matrices $\Omega_i$ and $\Pi_j$.
%
%%%, with all entries outside the submatrix with indices in $S_i \times T_j$ zeroed out; zeroing out these elements obviously does not change the weighted loss.
%
We then define the \emph{locally-denoised matrix}:
\begin{align}
\what \X^{\loc} = \sum_{i=1}^{I} \sum_{j=1}^{J} \Omega_i \what \X_{(i,j)}^{\loc} \Pi_j.
\end{align}
We summarize the localized denoising procedure in Algorithm \ref{alg-local}.

\begin{algorithm}[ht]
\caption{Localized denoising for unweighted loss}
\label{alg-local}
\begin{algorithmic}[1]
\item 
Input: $\Y$; pairwise orthogonal projections $\Omega_1,\dots,\Omega_I$, $\Pi_1,\dots,\Pi_J$ \begin{tabbing}
\hspace*{0cm}\=\kill
\> $\sum_{i=1}^{I} \Omega_i = \I_p$; $\sum_{j=1}^{J} \Pi_j = \I_n$
\end{tabbing}

\item
for $1 \le i \le I$, $1 \le j \le J$:
\begin{tabbing}
\hspace*{0cm}\=\kill
\>$\what \X_{(i,j)}^{\loc}$ is output of Algorithm \ref{alg:denoise} 
    with weights $\Omega_i$ and $\Pi_j$ \\
\>$\AMSE_{(i,j)}^{\loc}$ is estimated mean squared error
\end{tabbing}

\item
$\what \X^{\loc} = \sum_{i=1}^{I} \sum_{j=1}^{J} \Omega_i \what \X_{(i,j)}^{\loc} \Pi_j$

\item
$\AMSE^{\loc} = \sum_{i=1}^{I} \sum_{j=1}^{J} \AMSE_{(i,j)}^{\loc}$

\end{algorithmic}
\end{algorithm}

\subsection{Performance of localized denoising}
\label{sec-details-loc}

%%%We will denote by $\what \X^{\shr} = \sum_{k=1}^{r} t_k c_k \tilde c_k \hat \u_k \hat \v_k^T$ the matrix that is optimally shrunk with respect to ordinary Frobenius loss. We then have the following result:

The following results compare the behavior of the localized denoiser $\what \X^{\loc}$ to the optimal singular value shrinker $\what \X^{\shr}$.

%%%We recall from Section \ref{sec-het} that $\u_k$ is \emph{generic} with respect to $\Omega_i$ if $\u_k^T \Omega_i \u_k \sim \tr{\Omega_i} / p$, and is \emph{heterogeneous} with respect to $\Omega_i$ otherwise; similarly for $\v_k$ and $\Pi_j$.

%%%In both Theorem \ref{thm-compare1} and Theorem \ref{thm-compare2}, we suppose $\Omega_1,\dots,\Omega_I$ and $\Pi_1,\dots,\Pi_J$ are pairwise orthogonal projections satisfying $\I_p = \sum_{i=1}^{I} \Omega_i$ and $\I_n = \sum_{j=1}^{J} \Pi_j$. 

\begin{thm}
\label{thm-compare1}
\begin{math}
\|\what \X^{\loc} - \X\|_{\Fr}^2
\le \|\what \X^{\shr} - \X\|_{\Fr}^2
\end{math}
 almost surely as $p,n\to \infty$.
\end{thm}

\begin{thm}
\label{thm-compare2}
Suppose that either $\u_1,\dots,\u_r$ are weighted orthogonal with respect to all $\Omega_i$, or $\v_1,\dots,\v_r$ are weighted orthogonal with respect to all $\Pi_j$. Then almost surely as $p,n\to \infty$,
\begin{math}
\|\what \X^{\loc} - \X\|_{\Fr}^2
\le \|\what \X^{\shr} - \X\|_{\Fr}^2 - \xi,
\end{math}
where $\xi \ge 0$, and $\xi > 0$ if some $\u_k$ is heterogeneous with respect to some $\Omega_i$ or some $\v_k$ is heterogeneous with respect to some $\Pi_j$.
\end{thm}

In other words, unless all the $\u_k$ are generic with respect to all of the $\Omega_i$ \emph{and} all the $\v_k$ are generic with respect to all of the $\Pi_i$, localized denoising will outperform singular value shrinkage asymptotically. The proofs of Theorems \ref{thm-compare1} and \ref{thm-compare2} are found in Section \ref{appendix-compare1} and Section \ref{appendix-compare2}, respectively.

\begin{rmk}
The weighted orthogonality condition of Theorem \ref{thm-compare2} will hold if the columns of $\X$ are drawn iid from a sufficiently well-behaved distribution in $\R^p$.
%%%, as is often the case in applications.
\end{rmk}

\begin{rmk}
\label{rmk-projections}
To apply Theorem \ref{thm-compare2}, the user must select projection matrices $\Omega_i$ and $\Pi_j$ with respect to which the singular vectors of $\X$ are heterogeneous.
%%% In many applications, however, this knowledge may well be available. For instance, 
Datasets are often drawn from different experimental regimes in genetic microarray experiments \cite{johnson2007adjusting, leek2010tackling, scherer2009batch}, single-cell RNA processing \cite{shaham2017removal, tung2017batch}, and medical imaging \cite{kothari2011automatic}. In these settings, it is known a priori that signal components will likely be heterogeneous across the different subpopulations, and localized shrinkage is a natural tool.
%
%%%Furthermore, Figure \ref{fig-mit} illustrates that even when precise clustering is not used, localized denoising can still achieve substantial reduction in error compared to singular value shrinkage.
\end{rmk}

\begin{rmk}
\label{rmk-finite-sample}
Theorem \ref{thm-compare1} guarantees that even if the projection matrices $\Omega_i$ and $\Pi_j$ are not chosen judiciously (see Remark \ref{rmk-projections}), the asymptotic performance of localized denoising is never worse than singular value shrinkage. In practice, localized denoising requires estimating more parameters than does shrinkage, and if $I$ and $J$ are sizeable relative to $p$ and $n$ its performance might be worse due to finite sample fluctations in these parameter estimates, especially when the singular vectors of $\X$ do not exhibit strong heterogeneity with respect to the projections. For such an example, see Section \ref{sec-numerics-localized}, and specifically Remark \ref{rmk-local}. 
Though a detailed analysis of this topic is beyond the scope of the present work, in practice the user can compare these trade-offs via simulation to determine if localized denoising is appropriate for their problem size and the expected level of heterogeneity with respect to the projections.
\end{rmk}

\section{Applications of weighted denoising}
\label{sec-motivation}

In this section, we describe three applications of weighted loss functions: submatrix denoising, denoising with heteroscedastic noise, and denoising with missing data. In these problems, we estimate a low-rank matrix with respect to \emph{unweighted} Frobenius loss. However, an intermediate step of the estimation procedure requires the use of a \emph{weighted} loss function.

%%% but other aspects of the statistical model lead us to perform denoising with respect to a weighted loss function.

%%%That is, even when our ultimate goal is unweighted estimation, an intermediate step of the estimation procedure requires the use of a weighted loss function.

%

\subsection{Submatrix denoising}
\label{sec-splitting}

We suppose we observe a data matrix $\Y = \X + \G$, but our goal is to estimate only a $p_0$-by-$n_0$ submatrix of $\X$, where $p_0 / p \sim \mu$ and $n_0 /n \sim \nu$. Denoting by $\Omega \in \R^{p_0 \times p}$ the coordinate selection operator for the $p_0$ rows of the submatrix, and $\Pi \in \R^{n_0 \times n}$ the coordinate selection operator for the $n_0$ columns of the submatrix, we may write the target submatrix as $\X_0 = \Omega \X \Pi^T$.

One approach is to estimate the entire matrix $\X$ with respect to the weighted loss
\begin{math}
\L(\what \X, \X) = \| \Omega(\what \X - \X) \Pi^T\|_{\Fr}^2.
\end{math}
This loss only penalizes errors in the $p_0$ rows and $n_0$ columns of $\X_0$. If $\what \X$ denotes the optimal spectral denoiser minimizing $\L(\what \X, \X)$, we define our estimator $\what \X_0 = \Omega \what \X \Pi^T$. The method is summarized in Algorithm \ref{alg-submatrix}.

\begin{algorithm}[ht]
\caption{Submatrix denoising}
\label{alg-submatrix}
\begin{algorithmic}[1]
\item 
Input: $\Y$; projections $\Omega$, $\Pi$

\item
\begin{tabbing}
\hspace*{0cm}\=\kill
\>$\what \X$ is output of Algorithm \ref{alg:denoise} 
    with weights $\Omega$ and $\Pi$ \\
\>$\AMSE$ is estimated mean squared error
\end{tabbing}

\item
$\what \X_0 = \Omega \what \X \Pi^T$

\end{algorithmic}
\end{algorithm}

Another natural approach is to simply ignore the $p-p_0$ rows and $n-n_0$ columns outside of $\X_0$, and denoise $\X_0$ directly by optimal singular value shrinkage to the matrix $\Y_0 = \X_0 + \G_0$ (defining $\G_0 = \Omega \G \Pi^T$). We let $\what \X_0^{\shr}$ denote this estimator.

In the following result, we make the same assumptions on $\Omega$ and $\Pi$ from Section \ref{sec-model}. Note that $p_0 = \tr{\Omega^T \Omega}$, and $n_0 = \tr{\Pi^T \Pi}$.

\begin{prop}
\label{prop-merge}
Suppose $\u_1,\dots,\u_r$ are weighted orthogonal with respect to $\Omega^T \Omega$, and $\v_1,\dots,\v_r$ are weighted orthogonal with respect to $\Pi^T \Pi$. Suppose
\begin{math}
\alpha_k < \sqrt{\mu}
\end{math}
and
\begin{math}
\beta_k < \sqrt{\nu},
\end{math}
for $1 \le k \le r$.
Then 
\begin{math}
\|\what \X_0 - \X_0\|_{\Fr}^2
< \|\what \X_0^{\shr} - \X_0\|_{\Fr}^2,
\end{math}
where the strict inequality holds almost surely in the limit $p,n \to \infty$.
\end{prop}

The proof of Proposition \ref{prop-merge} is found in Section \ref{appendix-merge}.

\begin{rmk}
If $\u_k$ and $\v_k$ are generic with respect to $\Omega^T \Omega$ and $\Pi^T \Pi$, respectively, then $\alpha_k = \mu$ and $\beta_k = \nu$. Proposition \ref{prop-merge} requires the much weaker condition that $\alpha_k \le \sqrt{\mu}$ and $\beta_k \le \sqrt{\nu}$ (note that $\mu \le \sqrt{\mu}$ and $\nu \le \sqrt{\nu}$). Informally, even if the fraction of the signal's energy contained in $\X_0$ is disproportionately large, it still pays to denoise $\X_0$ using the entire observed matrix $\Y$, rather than the submatrix $\Y_0$ alone.
%%%In other words, the $p-p_0$ rows and $n - n_0$ columns of $\Y$ lying outside $\Y_0$ serve to boost the signal in the submatrix $\Y_0$.
\end{rmk}

%%%\begin{rmk}
It will follow from the proof of Proposition \ref{prop-merge} that if $\alpha_k < \sqrt{\mu}$ and $\beta_k < \sqrt{\nu}$, then the singular vectors of $\X_0$ are better approximated by computing the singular vectors of $\Y$ and projecting onto the images of $\Omega$ and $\Pi$, respectively, rather than computing the singular vectors of the submatrix $\Y_0$ itself. More precisely, we will show that 
\begin{math}
\u_k^0 = \frac{\Omega \u_k}{\|\Omega \u_k\|}
\end{math}
and
\begin{math}
\v_k^0 = \frac{\Pi \v_k}{\|\Pi \v_k\|}
\end{math}
are the singular vectors of $\X_0$, and that the vectors
\begin{math}
\hat \u_k^\omega = \frac{\Omega \hat \u_k}{\|\Omega \hat \u_k\|}
\end{math}
and
\begin{math}
\hat \v_k^\omega = \frac{\Pi \hat \v_k}{\|\Pi \hat \v_k\|},
\end{math}
are better correlated with $\u_k^0$ and $\v_k^0$, respectively, then are the singular vectors of $\Y_0$.
%%%\end{rmk}

%

\subsection{Doubly-heteroscedastic noise}
\label{sec-whitening}

We consider estimating a low-rank matrix $\X$ from an observed matrix $\Y = \X + \N$, where $\N$ is a noise matrix of the form $\N = \S^{1/2} \G \T^{1/2}$, $\G$ has iid entries with distribution $N(0,1/n)$, and $\S \in \R^{p \times p}$ and $\T \in \R^{n \times n}$ are positive-definite matrices. We assume the eigenvalues of $\S = \S_p$ and $\T = \T_n$ remain in an interval $[a,b]$ for all $p$ and $n$, where $a>0$ and $b < \infty$ are fixed independently of $p$ and $n$. We refer to the matrix $\N$ as \emph{doubly-heteroscedastic} noise.

\begin{rmk}
This noise model generalizes two previous models of heteroscedastic noise in the context of principal component analysis \cite{hong2016towards, hong2018asymptotic, hong2018optimally, zhang2018heteroskedastic, leeb2019optimal}. In both, the matrix $\X$ consists of random, iid signal vectors $X_1,\dots,X_n$ of the form
\begin{math}
X_j = \sum_{k=1}^{r} \ell_k^{1/2} z_{jk} \u_k,
\end{math}
where $\ell_1 > \dots > \ell_r > 0$, $\u_1,\dots,\u_r$ are orthonormal vectors, and the $z_{jk}$ are iid random variables with variance $1$ and mean $0$.
%%%The $\u_1,\dots,\u_r$ are the \emph{principal components} of the signal distribution.
The model from \cite{zhang2018heteroskedastic} and \cite{leeb2019optimal} takes $\T = \I_n$, in which case the observations are of the form $Y_j = X_j +  \S^{1/2} G_j$, where $G_j \sim N(\mathbf{0},\I_p)$.
%%%In this model, each noise vector $\S^{1/2} G_j$ has the same distribution, whereas the noise level may vary across each of the $p$ coordinates.
By contrast, the papers \cite{hong2016towards, hong2018asymptotic, hong2018optimally} take $\S = \I_p$, in which case the observations are of the form $Y_j = X_j +  b_j^{1/2} G_j$.
%%%Here, each noise vector $b_j^{1/2} G_j$ is isotropic, but the noise strength varies across the observations.
The doubly-heteroscedastic noise model takes $Y_j = X_j + b_j^{1/2} \S^{1/2} G_j$, which generalizes both these models.

%%%We note too that random matrices with a rank $1$ variance profile have also been studied in \cite{hachem2007deterministic, hachem2008CLT, leeb2020rapid}.
\end{rmk}

We consider the following three-step procedure. First, we \emph{whiten} the noise, replacing $\Y$ by $\wtilde \Y$ defined by
\begin{math}
\wtilde \Y = \S^{-1/2} \Y \T^{-1/2}.
\end{math}
We may write $\wtilde \Y = \wtilde \X + \G$, where $\wtilde \X = \S^{-1/2} \X \T^{-1/2}$ and $\G$ has iid $N(0,1/n)$ entries. Next, we apply a denoiser to $\wtilde \Y$ to estimate $\wtilde \X$; we denote this by $\psi(\wtilde \Y)$, for $\psi$ tailored to removing \emph{white} noise. Finally, we \emph{unwhiten} $\psi(\wtilde \Y)$ to obtain our final estimate $\what \X = \S^{1/2} \psi(\wtilde \Y) \T^{1/2} $ of $\X$.

The Frobenius loss between $\what \X$ and $\X$ may be written as follows:
\begin{align}
\|\what \X - \X   \|_\Fr^2 
= \| \S^{1/2} \psi(\wtilde \Y) \T^{1/2} - \S^{1/2}\wtilde\X\T^{1/2}  \|_\Fr^2
= \| \S^{1/2} (\psi(\wtilde \Y) - \wtilde\X ) \T^{1/2}  \|_\Fr^2,
\end{align}
which is a \emph{weighted} loss between $\psi(\wtilde \Y)$ and $\wtilde \X$, with weights $\S^{1/2}$ and $\T^{1/2}$. The denoiser $\psi(\wtilde \Y)$ should be chosen to minimize this weighted Frobenius loss. The procedure is summarized in Algorithm \ref{alg-hetero}, where $\psi$ is taken to be optimal spectral denoising,

\begin{algorithm}[ht]
\caption{Matrix denoising with doubly-heteroscedastic noise}
\label{alg-hetero}
\begin{algorithmic}[1]
\item 
Input: $\Y$; positive-definite $\S$, $\T$

\item
$\wtilde \Y = \S^{-1/2} \Y \T^{-1/2}$
\item
\begin{tabbing}
\hspace*{0cm}\=\kill
\>$\psi(\wtilde \Y)$ is output of Algorithm \ref{alg:denoise} 
    with weights $\S^{1/2}$ and $\T^{1/2}$ \\
\>$\AMSE$ is estimated mean squared error
\end{tabbing}

\item
$\what \X = \S^{1/2} \psi(\wtilde \Y) \T^{1/2}$

\end{algorithmic}
\end{algorithm}

\begin{rmk}
The procedure of whitening, denoising, and unwhitening has been employed in recent papers on the spiked model; see, for instance, \cite{liu2016epca, leeb2019optimal, dobriban2017optimal}. In particular, \cite{leeb2019optimal} shows several advantages of working with the whitened matrix when the noise is one-sided, such as improved estimation of the singular vectors of $\X$.
%%%In Section \ref{sec-snr}, namely Proposition \ref{prop-snr}, we show that whitening doubly-heteroscedastic noise increases the signal-to-noise ratio of the observed matrix.
By contrast, the paper \cite{hong2018optimally} shows that whitening is suboptimal in certain settings.
\end{rmk}

\subsubsection{Estimating $\S$ and $\T$}
The signal/noise decomposition $\Y = \X + \N$ is obviously not well-defined unless the user possesses some additional knowledge about the noise matrix $\N = \S^{1/2} \G \T^{1/2}$.  While a detailed treatment of this problem is outside the scope of this paper, we observe that in the large $p$, large $n$ asymptotic limit, the matrices $\S = \S_p$ and $\T = \T_n$ may be replaced by estimators $\what \S = \what \S_p $ and $\what \T = \what \T_n$ consistent in operator norm; that is, almost surely
\begin{align}
\label{consistency-operator}
\lim_{p \to \infty} \| \S_p - \what \S_p\|_{\op} = 
\lim_{n \to \infty} \| \T_n - \what \T_n\|_{\op} = 0.
\end{align}

\begin{rmk}
The matrices $\S$ and $\T$ may be replaced by, respectively, $\theta\S$ and $\T/\theta$ for any $\theta > 0$. Without loss of generality, we may therefore assume that $\tr{\T}/n = 1$.
\end{rmk}

%%%Replacing $\S$ and $\T$ with $\what \S$ and $\what \T$, respectively, in Algorithm \ref{alg-hetero}, the output matrix $\what \X$ w

The next result describes a simple method for estimating $\S$ and $\T$ consistently in operator norm when both are diagonal and the singular vectors of $\X$ are delocalized.

\begin{prop}
\label{prop-estimating}
Suppose  $\max_{1 \le k \le r}\|\u_k\|_\infty \|\v_k\|_\infty = o(n^{-1/2})$, $\S = \diag(a_1,\dots,a_p)$ and $\T = \diag(b_1,\dots,b_n)$. For $1 \le i \le p$ and $1 \le j \le n$, define the estimators
\begin{align}
\label{ahat}
\hat a_i = \sum_{j=1}^{n} Y_{ij}^2,
\quad
\hat b_j = \frac{\sum_{i=1}^{p} Y_{ij}^2}{\frac{1}{n}\sum_{i=1}^{p} \hat a_i },
\end{align}
and
and let $\what \S = \what \S_p = \diag(\hat a_1,\dots,\hat a_p)$ and $\what \T = \what \T_n = \diag(\hat b_1,\dots,\hat b_p)$. Then $\what \S$ and $\what \T$ are consistent estimators of $\S$ and $\T$, respectively; that is, \eqref{consistency-operator} holds almost surely.
\end{prop}

The proof of Proposition \ref{prop-estimating} may be found in Section \ref{appendix-estimating}.

\begin{rmk}
The values $\hat a_i$ in \eqref{ahat} are the sample standard deviations of the rows of $\sqrt{n}\Y$. Normalizing $\Y$ by $\what \S^{1/2}$ is then an instance of \emph{standardization} of the rows, a commonly used method in principal component analysis \cite{jolliffe2002principal}.
\end{rmk}

\begin{rmk}
The estimates $\hat a_i$ and $\hat b_j$ capture the variation of both the noise and the signal. Proposition \ref{prop-estimating} states that if the signal is delocalized, then in the large $p$, large $n$ limit its contribution becomes negligible. However, for finite $p$ and $n$, $\what \S$ and $\what \T$ will still see the effects of the signal, and may not be good estimates of $\S$ and $\T$. The experiment in Section \ref{sec-numerics-whitening} compares the use of the true $\S$ and $\T$ to their estimates.
\end{rmk}

\subsubsection{Whitening increases the SNR for generic signal matrices}
\label{sec-snr}
We show that the whitening transformation increases a natural signal-to-noise ratio of the observed matrix. We will assume throughout this section that the $\u_k$ (respectively, $\v_k$) are generic with respect to $\S$ (respectively, $\T$), and that they satisfy the pairwise orthogonality condition with respect to $\S$ (respectively $\T$). Writing the SVD of $\X$ as
\begin{math}
\X = \sum_{k=1}^{r} t_k \u_k \v_k^T,
\end{math}
we define the signal-to-noise ratio (SNR) for component $k$ of $\X$:
\begin{align}
\snr_k = \frac{ t_k^2 }{\| \N \|_{\op}^2},
\end{align}
which is the ratio of the squared operator norm of the component $t_k \u_k \v_k^T$ of $\X$ and the squared operator norm of the noise.

%%% $\snr_k$ measures the strength of the $k^{th}$ component relative to the strength of the noise matrix $\N$.

Whitening turns $\Y$ into $\wtilde \Y = \wtilde \X + \G$, with
\begin{math}
\wtilde \X = \sum_{k=1}^{r} t_k (\S^{-1/2} \u_k) (\T^{-1/2}\v_k)^T
= \sum_{k=1}^{r} \tilde t_k \tilde \u_k \tilde \v_k^T,
\end{math}
where $\tilde t_k = t_k \| \S^{-1/2} \u_k \| \| \T^{-1/2} \v_k \|$, $\tilde \u_k = \S^{-1/2} \u_k / \|\S^{-1/2} \u_k\|$, and $\tilde \v_k = \T^{-1/2} \v_k / \|\T^{-1/2} \v_k\|$.  The  SNR after whitening is then:
\begin{align}
\wtilde \snr_k = \frac{ \tilde t_k^2 }{\| \G \|_{\op}^2},
\end{align}

Define
\begin{align}
\tau = \left(\frac{1}{p}\tr{\S}\right) \left(\frac{1}{p}\tr{\S^{-1}}\right)
    \left(\frac{1}{n}\tr{\T}\right) \left(\frac{1}{n}\tr{\T^{-1}}\right).
\end{align}
Note that from Jensen's inequality, $\tau \ge 1$, and $\tau > 1$ if either $\S$ or $\T$ is not a multiple of the identity. The following result extends an analogous finding from \cite{leeb2019optimal}:
\begin{prop}
\label{prop-snr}
Suppose $\u_1,\dots,\u_r$ are generic and weighted orthogonal with respect to $\S$, and $\v_1,\dots,\v_r$ are generic and weighted orthogonal  with respect to $\T$. Then
\begin{math}
\wtilde \snr_k \ge \tau \snr_k,
\end{math}
$1 \le k \le r$,
almost surely as $p,n\to \infty$. In particular, $\wtilde \snr_k$ is larger than $\snr_k$ if either $\S$ or $\T$ is not a multiple of the identity.
\end{prop}

In other words, the SNR increases after whitening the noise by at least a factor of $\tau$; energy is transferred from the noise component to the signal component. The proof of Proposition \ref{prop-snr}, which extends an analogous result in \cite{leeb2019optimal}, is in Section \ref{appendix-snr}.

\subsection{Matrices with missing/unobserved values}
\label{sec-missing}

We consider the setting where $\X$ is a low-rank target matrix we wish to recover and $\G$ is a matrix of iid Gaussian $N(0,1)$ entries, but rather than observe $\X + \G$, we observe only some subset of the entries. The problem of estimating a matrix from a subset of its entries is known as \emph{matrix completion} \cite{recht2011simpler, keshavan2010matrix, jain2013, candes2010power, chen2015completing, candes2009exact, candes2010power, klopp2014noisy, keshavan2010regularization, candes2010noise, dobriban2017optimal, koltchinskii2011nuclear, negahban2011estimation, srebro2010collaborative}.

In this section, we will adopt a heterogeneous, rank 1 sampling model, as in \cite{chen2015completing}. We suppose that the rows and columns are sampled independently, with row $i$ sampled with probability $q_i^r$, and column $j$ sampled with probability $q_j^c$; that is, entry $(i,j)$ of $\X+\G$ is sampled with probability $q_i^r q_j^c$. We observe the vector  $\y = \F(\X+\G)$ of $M$ sampled entries, where $\F: \R^{p \times n} \to \R^M$ is the subsampling operator.
%%% $M$ being the number of sampled entries.

Following the approach from \cite{dobriban2017optimal}, we consider the backprojected matrix $\Y = \F^* (\y) \in \R^{p \times n}$, in which the unobserved entries are replaced by $0$'s. We write $\Y = \F^* (\F (\X)) + \F^* (\F (\G))$. We show that asymptotically, $\F^* (\F (\X))$ behaves like the matrix $\P \X \Q$. More precisely, we have the following result:
\begin{prop}
\label{prop-missing}
Suppose $\max_{1 \le k \le r}\|\u_k\|_\infty \|\v_k\|_\infty = o(n^{-1/2})$. Then in the limit $p / n \to \gamma$, 
\begin{math}
\| \F^* (\F (\X)) - \P \X \Q \|_{\op} \to 0
\end{math}
almost surely.
\end{prop}
The proof of Proposition \ref{prop-missing} may be found in Section \ref{appendix-missing}. It is a straightforward generalization of the analogous one-sided result in \cite{dobriban2017optimal}.

Let $\N = \F^*(\F(\G))$. Writing $\N = (N_{ij})$, we have $N_{ij} = \delta_{ij} G_{ij}$, where $\delta_{ij}$ is $1$ if entry $(i,j)$ is sampled, and $0$ otherwise. Then $N_{ij}$ has variance $q_i^r q_j^c$. Consequently, we  can whiten the noise by applying $\P^{-1/2}$ and $\Q^{-1/2}$; Proposition \ref{prop-snr} suggests this will improve estimation of the matrix. To that end, we define
\begin{math}
\wtilde \Y = \P^{-1/2} \Y \Q^{-1/2} = \wtilde \X + \wtilde \G,
\end{math}
where $\wtilde \X = \P^{-1/2}\F^*(\F(\X)) \Q^{-1/2}$, and $\wtilde \G = \P^{-1/2} \N \Q^{-1/2}$. Then $\wtilde \G$ is a random matrix where each entry has mean zero and variance $1$.

From Proposition \ref{prop-missing}, asymptotically the matrix $\wtilde \X$ behaves like $\P^{1/2} \X \Q^{1/2}$, and so denoising $\wtilde \Y$  estimates $\psi(\wtilde \Y)$ of $\P^{1/2} \X \Q^{1/2}$. To estimate $\X$ we should perform denoising to $\wtilde \Y$ with respect to the weighted loss function 
\begin{math}
\L(\what \X,\X) = \|\P^{-1/2} (\what \X - \X) \Q^{-1/2}\|_{\Fr}^2,
\end{math}
with weight matrices $\P^{-1/2}$ and $\Q^{-1/2}$. We then apply $\P^{-1/2}$ and $\Q^{-1/2}$ to the resulting matrix, to obtain an estimator of $\X$ itself. The method is summarized in Algorithm \ref{alg-missing}, where $\psi$ is the optimal spectral denoiser.

\begin{algorithm}[ht]
\caption{Matrix denoising with missing data}
\label{alg-missing}
\begin{algorithmic}[1]
\item 
Input: Samples $\y$; sampling operator $\F$; sampling matrices $\P$, $\Q$

\item
Normalize, backproject observations $\wtilde \Y = \P^{-1/2} \F^*(\y) \Q^{-1/2}$

\item
\begin{tabbing}
\hspace*{0cm}\=\kill
\>$\psi(\wtilde \Y)$ is output of Algorithm \ref{alg:denoise} 
    with matrix $\wtilde \Y$, weights $\P^{-1/2}$ and $\Q^{-1/2}$ \\
\>$\AMSE$ is estimated mean squared error
\end{tabbing}

\item
$\what \X = \P^{-1/2} \psi(\wtilde \Y) \Q^{-1/2}$

\end{algorithmic}
\end{algorithm}

%

%
%    figure from experiments/expt_local2
%
\begin{figure}[h]
\centering
\includegraphics[scale=0.4]{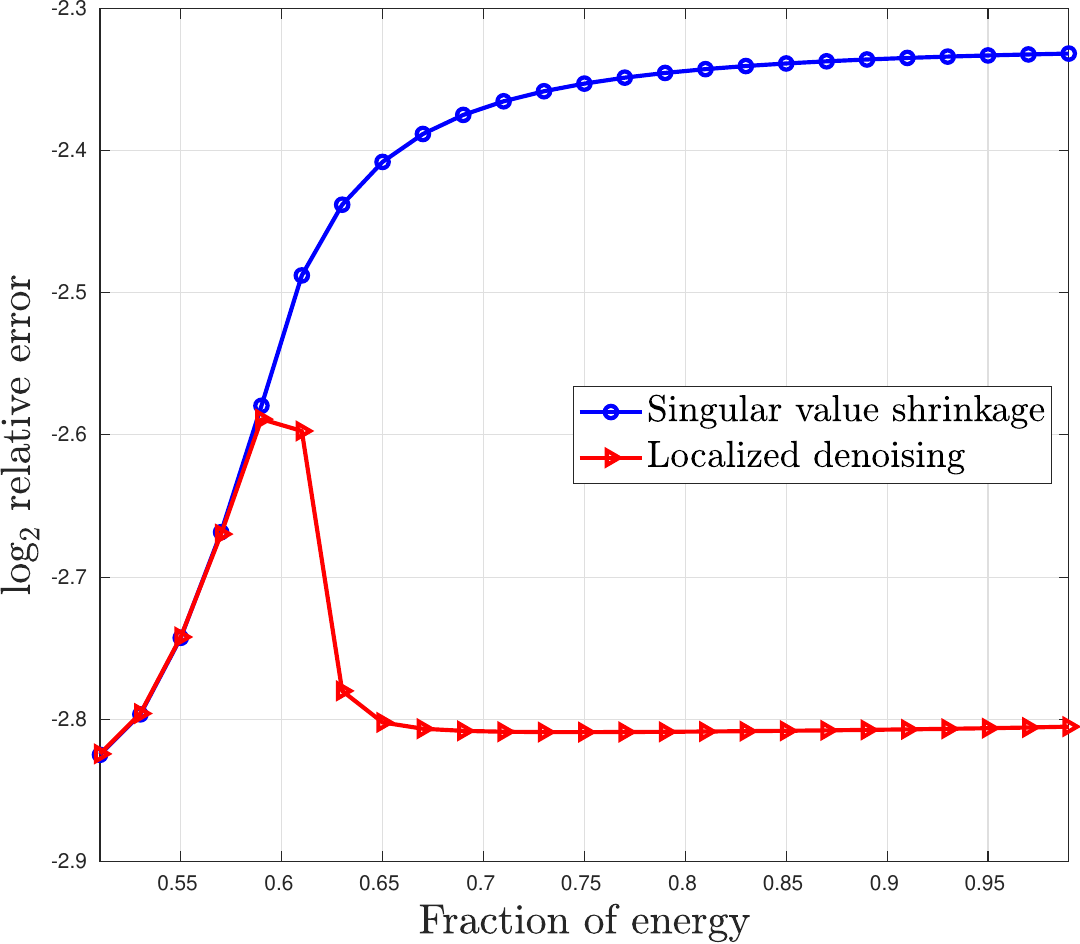}
\caption{
Localized denoising versus singular value shrinkage for the checkerboard matrix, shown in Figure \ref{fig-checkers}; see Section \ref{sec-numerics-localized} for simulation details. The $x$-axis is parametrized by the fraction of signal energy contained in the light squares. Localized denoising outperforms shrinkage when the fraction is large enough that the rank $2$ block structure is detectable.
}
\label{fig-local}
\end{figure}

%
%    figure from experiments/checkers_denoised2
%
\begin{figure}[h]
\centering
\includegraphics[scale=0.4]{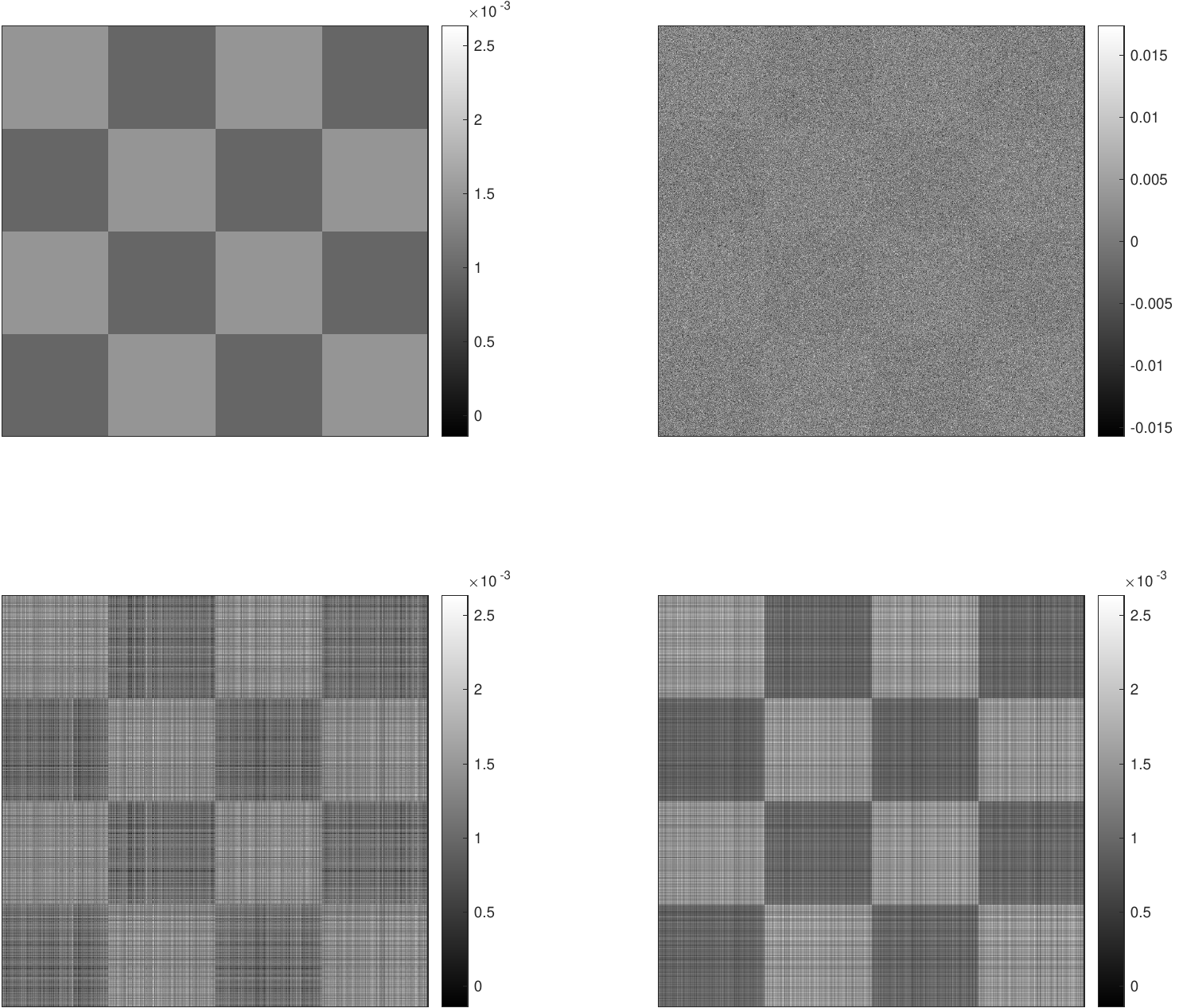}
\caption{
Localized denoising versus singular value shrinkage; see Section \ref{sec-numerics-localized} for simulation details. Upper left: the rank $2$ signal matrix; the fraction of signal energy in the light squares is $f = 0.7$. Upper right: the observed noisy matrix. Lower left: the matrix denoised by optimal singular value shrinkage \cite{gavish2017optimal, shabalin2013reconstruction}. Lower right: the matrix denoised by localized denoising (Algorithm \ref{alg-local}). The relative error of singular value shrinkage is approximately $1.92 \times 10^{-1}$, whereas the relative error of localized denoising is approximately $1.40 \times 10^{-1}$.
}
\label{fig-checkers}
\end{figure}

\section{Numerical results}
\label{sec-numerics}

In this section, we report on numerical simulations demonstrating the performance of the algorithms from this paper.

%%%In Section \ref{sec-numerics-localized}, we compare localized denoising to global singular value shrinkage. In Section \ref{sec-numerics-submatrix}, we report on spectral denoising for submatrix denoising, as described in Section \ref{sec-splitting}. In Section \ref{sec-numerics-whitening}, we report on spectral denoising for doubly-heteroscedastic noise, as described in Section \ref{sec-whitening}. In Section \ref{sec-numerics-missing}, we report on spectral denoising method for missing data, as described in Section \ref{sec-missing}. In Section \ref{sec-nongaussian}, we examine the accuracy of the formulas from Theorem \ref{thm-asymptotics} for non-Gaussian noise matrices $\G$. In Section \ref{sec-rank}, we explore the question of estimating the rank $r$ of $\X$.

%

%
%    figure from experiments/expt_submatrix2
%
\begin{figure}[h]
\centering
\includegraphics[scale=0.4]{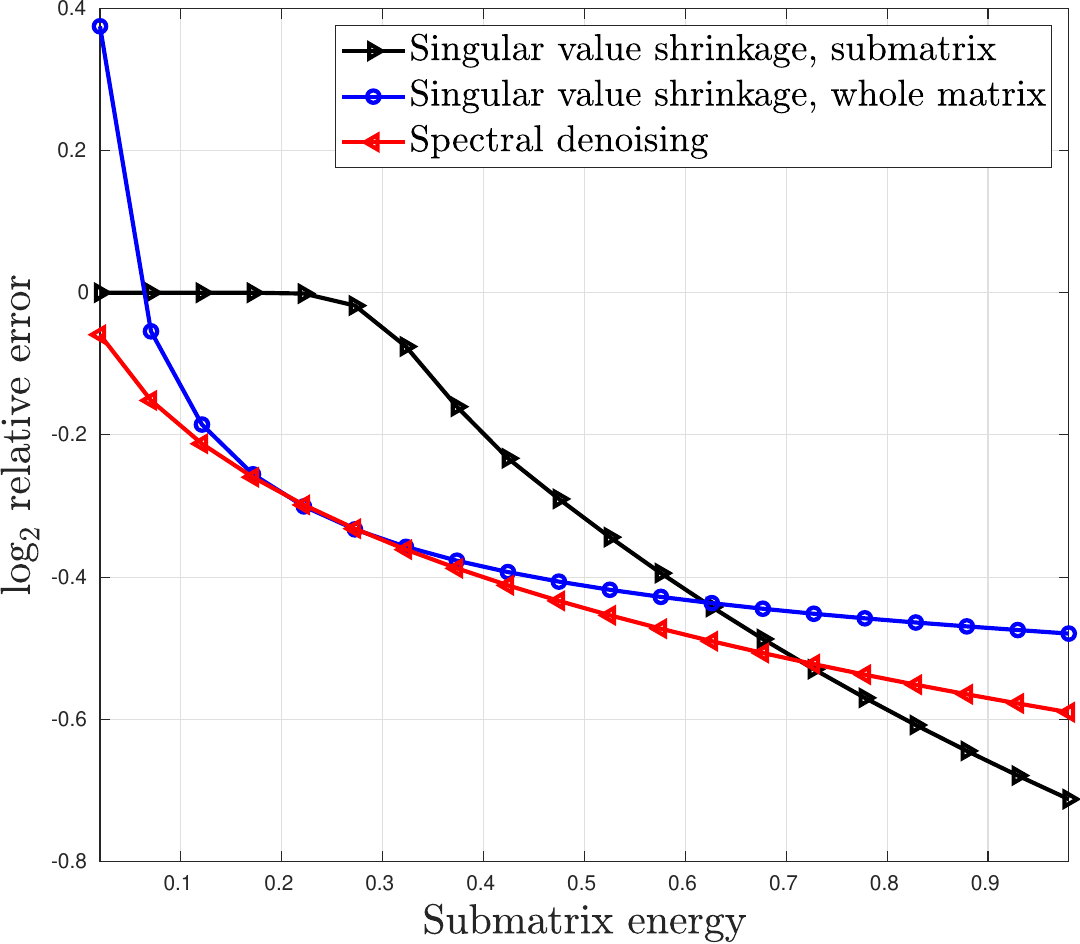}
\caption{
Submatrix denoising; see Section \ref{sec-numerics-submatrix} for simulation details. The left and right singular vectors of $\X$ each contain $\sqrt{f} \times 100\%$ of their energy in the submatrix coordinates, and the $x$-axis is parametrized by the fraction of energy $f=\|\X_0\|_{\Fr}^2 / \||\X\|_{\Fr}^2$ in the submatrix $\X_0$. The $\log_2$ relative error of singular value shrinkage on the submatrix plateaus to $0$ when $f$ is small, since the signal in the submatrix alone is too weak to be detected. Singular value shrinkage on the submatrix outperforms spectral denoising (Algorithm \ref{alg-submatrix}) when $f$ is very large, but is otherwise inferior; this is the behavior expected from Proposition \ref{prop-merge}.
}
\label{fig-submatrix}
\end{figure}

\subsection{Localized denoising}
\label{sec-numerics-localized}

We evaluate the performance of localized denoising (Algorithm \ref{alg-local}). 
%
%%%For comparison, we also apply optimal singular value shrinkage \cite{gavish2017optimal, shabalin2013reconstruction} to the entire matrix $\X$. Theorem \ref{thm-compare2} predicts that localized denoising will outperform shrinkage when the matrix $\X$ is heterogeneous, i.e.\ when the energy in the singular vectors is not uniformly distributed across the coordinates.
%
We generate a ``checkerboard'' signal matrix $\X$ of size $p$-by-$n$, $p=n=800$, shown in the top left panel of Figure \ref{fig-checkers}. Each light square has the same value, as does each dark square. For a specified number $f \in [1/2,1]$,  the total energy of the light squares is $f \times 100\%$ of the total energy of $\X$. The Frobenius norm of $\X$ is normalized to be $1$. Whenever $f > 1/2$,  $\X$ is rank $2$; when $f=1/2$, $\X$ has constant value and is rank $1$. We add a matrix $\G$ of Gaussian noise, whose entries have standard deviation $1/(10 \sqrt{n})$.

We estimate $\X$ from $\Y$ using two methods: singular value shrinkage \cite{gavish2017optimal, shabalin2013reconstruction} and localized denoising. Localized denoising is applied with row projection matrices $\Omega_i$, $i=1,2,3,4$, that project onto equispaced blocks of rows, and column projection matrices $\Pi_i$, $i=1,2,3,4$, that project onto equispaced blocks of columns. For each $f$, the experiment is repeated $50$ times; the $\log_2$ mean errors are plotted in Figure \ref{fig-local}.

As $f$ increases, localized denoising outperforms singular value shrinkage more dramatically. This is because localized denoising uses a priori knowledge of $\X$'s block structure, which becomes more pronounced as $f$ grows. Figure \ref{fig-checkers} shows an example of images of the true matrix $\X$, the noisy matrix $\Y$, and the two denoised matrices, when $f = 0.7$. In this example, the relative error $\|\what \X^{\loc} - \X\|_{\Fr} / \|\X\|_{\Fr}$ of localized denoising is approximately $1.40 \times 10^{-1}$, whereas the shrinkage error is $1.92 \times 10^{-1}$.

%%%Additionally, the checkerboard image is much sharper for localized denoising.

\begin{rmk}
\label{rmk-local}
The error curves in Figure \ref{fig-local} both appear nearly identical when $f \lesssim 0.6$, after which localized denoising begins to outperform singular value shrinkage. This is because for small values of $f$  the smallest singular value of $\X$ is not detectable, and so both methods treat the matrix as a constant, rank $1$ matrix. Though not apparent from the plot, when $f \le 0.55$ the performance of singular value shrinkage is slightly better than localized denoising, due to finite sample fluctuations (see Remark \ref{rmk-finite-sample}). For example, when $f=0.51$, the mean relative error of localized denoising is approximately $1.4118 \times 10^{-1}$, while that of shrinkage is approximately $1.4112 \times 10^{-1}$.
\end{rmk}

%
%    figure from experiments/expt_heter3
%
\begin{figure}[h]
\centering
\includegraphics[scale=0.4]{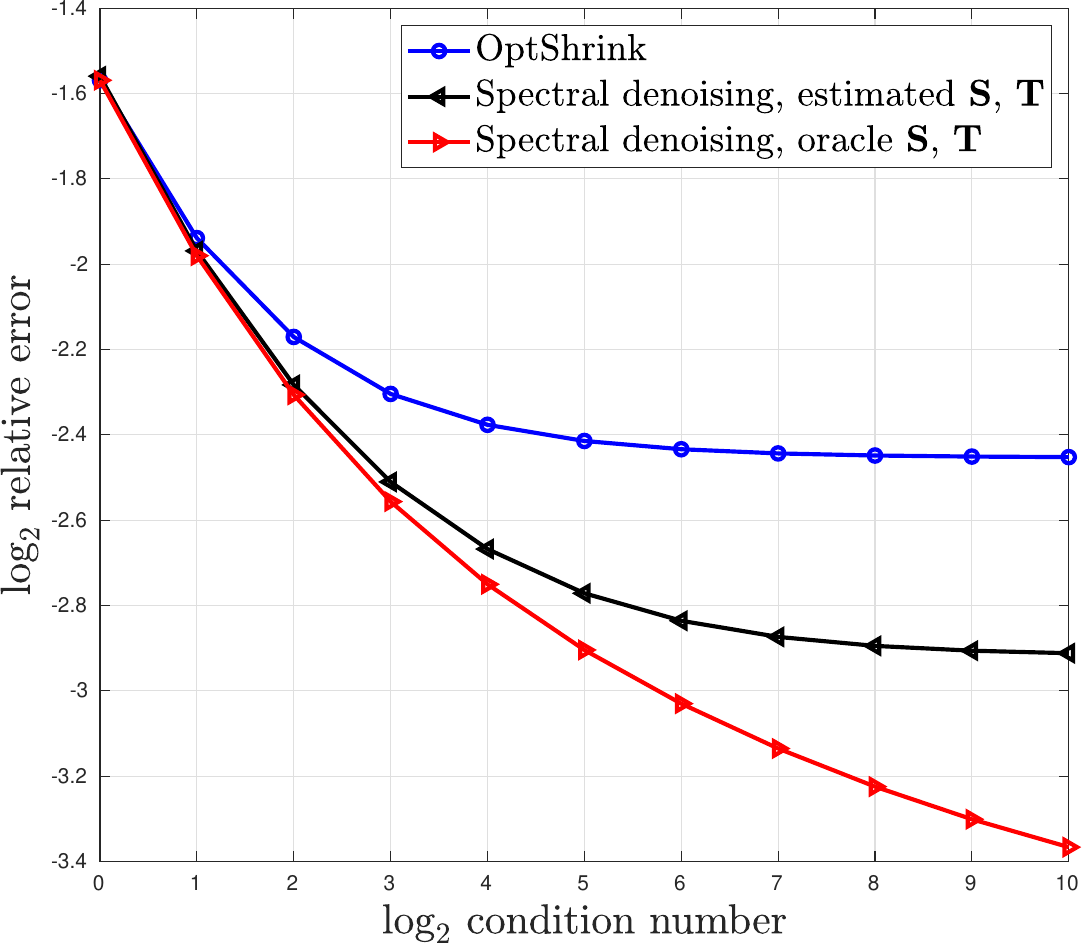}
\caption{
Denoising a matrix with doubly-heteroscedastic noise; see Section \ref{sec-numerics-whitening} for simulation details. The $x$-axis is parameterized by $\log_2$ of the condition number $\kappa$ of $\S$ and $\T$. Proposition \ref{prop-snr} suggests that noise whitening will enhance performance, and increasingly so as the condition number $\kappa$ grows, as is the case comparing spectral denoising with oracle whitening (Algorithm \ref{alg-hetero}) and OptShrink \cite{nadakuditi2014optshrink}. Interestingly, this appears to still hold even when $\S$ and $\T$ are estimated using the procedure from Proposition \ref{prop-estimating}.
}
\label{fig-heter}
\end{figure}

\subsection{Submatrix denoising}
\label{sec-numerics-submatrix}

We evaluate the performance of spectral denoising for estimating a submatrix $\X_0$ contained within a larger matrix $\X$ (Algorithm \ref{alg-submatrix}). We generate a rank 1 signal matrix $\X$ of size $p$-by-$n$, $p=500$, $n=1000$, with singular values $\gamma^{1/4} + 1/2$, where $\gamma=1/4$. For a specified $f\in (0,1)$, the left singular vector $\u=\u_1$ of $\X$ is piecewise constant on the two sets of coordinates $\{1,\dots,p/2\}$ and $\{p/2+1,\dots,p\}$; the values are such that the energy of $\u$ on coordinates $\{1,\dots,p/2\}$ is equal to $\sqrt{f}$. Similarly, the right singular vector $\v=\v_1$ of $\X$ is piecewise constant on the two sets of coordinates $\{1,\dots,n/2\}$ and $\{n/2+1,\dots,n\}$, with values such that the energy of $\v$ on $\{1,\dots,n/2\}$ is also equal to $\sqrt{f}$. Denoting by $\X_0$ the $p/2$-by-$n/2$ upper-left submatrix of $\X$, $f = \|\X_0\|_{\Fr}^2 / \|\X\|_{\Fr}^2$.

%%%; that is, $f$ is the fraction of $\X$'s energy contained in $\X_0$.

The noise matrix has Gaussian entries with variance $1/n$. We denoise the submatrix $\X_0$ using Algorithm \ref{alg-submatrix}; optimal singular value shrinkage \cite{gavish2017optimal, shabalin2013reconstruction} on $\X_0$ alone (``submatrix shrinkage''); and optimal singular value shrinkage on  $\X$ followed by projection onto the rows and columns of $\X_0$ (``global shrinkage''). For each $f$, the experiment is repeated for $50$ draws. Figure \ref{fig-submatrix} plots the $\log_2$ mean relative errors.

Optimal spectral denoising outperforms global shrinkage for all $f$, since singular value shrinkage is an instance of spectral denoising and hence will not do better than the optimal spectral denoiser. Optimal spectral denoising and global shrinkage perform nearly identically when $f \approx 1/4$, since in this regime the singular vectors of $\X$ are constant, and hence generic with respect to the weight matrices.

For small $f$, the relative error of global shrinkage exceeds $1$, since the submatrix's energy is very small compared to the rest of the matrix. By contrast, optimal spectral denoising with weights $\Omega$ and $\Pi$ highlights the rows and columns in $\X_0$.

Optimal spectral denoising outperforms submatrix shrinkage except when $f$ is  close to $1$. This is consistent with Proposition \ref{prop-merge}, which states that unless the energy of $\X$'s singular vectors are highly concentrated in the rows and columns of $\X_0$, optimal spectral denoising will outperform singular value shrinkage on the submatrix.

Finally, optimal singular value shrinkage on the submatrix has relative error $1$ when $f$ is small. This is because singular value shrinkage on the submatrix only computes the SVD of $\Y_0$, not $\Y$; when the energy in the submatrix $\X_0$ is too weak (i.e.\ $f$ is too small), no signal will be detected in the submatrix $\Y_0$ alone. By contrast, the singular values of the full matrix $\Y$ always reveal the presence of signal.

%
%    figure from experiments/expt_missing
%
\begin{figure}[h]
\centering
\includegraphics[scale=0.4]{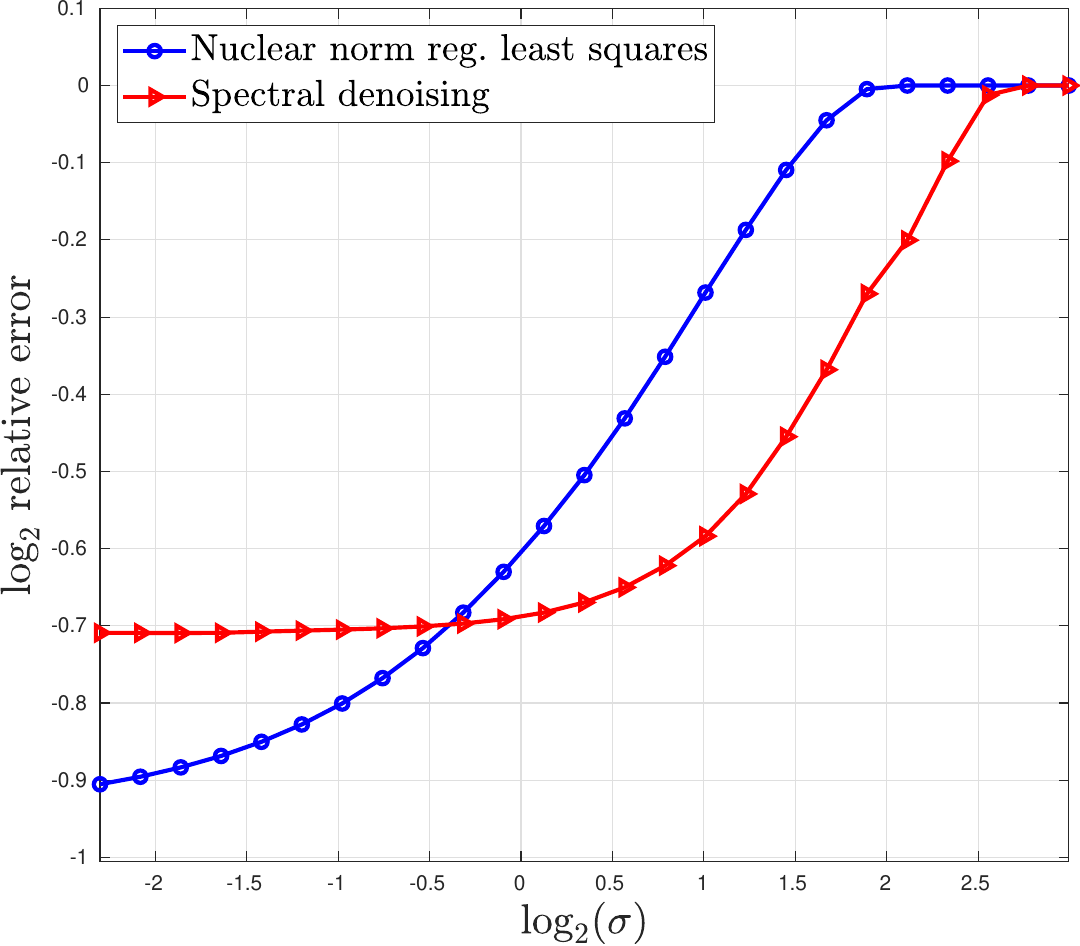}
\caption{
Denoising with missing data; see Section \ref{sec-numerics-missing} for simulation details. The $x$-axis is parameterized by the $\log_2$ noise level. The $\log_2$ relative errors of both methods plateau to $0$ at large $\sigma$, because the signal is undetectable in this regime. Spectral denoising for missing data (Algorithm \ref{alg-missing}) does better when $\sigma$ is large, but underperforms when $\sigma$ is small.
}
\label{fig-missing}
\end{figure}

\subsection{Doubly-heteroscedastic noise}
\label{sec-numerics-whitening}

We examine the performance of Algorithm \ref{alg-hetero}. We generate a $p$-by-$n$  signal matrix $\X$, $p=1000$, $n=2000$, of rank $r=5$, with singular values $t^* + 1/2 + k$, $k=0,1,2,3,4$, where $t^*$ is the smallest singular detectable value of $\X$, evaluated using the method in \cite{leeb2020rapid}. Both the left and right singular vectors of $\X$ are random orthonormal vectors in $\R^p$ and $\R^n$, respectively.

For specified $\kappa \ge 1$, we generate row and column diagonal covariance matrices $\S$ and $\T$, each with eigenvalues linearly spaced between $1/\kappa$ and $1$. The noise matrix is $\S^{1/2} \G \T^{1/2}$, where $\G$ has iid Gaussian entries with variance $1/n$. We apply three denoising schemes: Algorithm \ref{alg-hetero} with the true $\S$ and $\T$; Algorithm \ref{alg-hetero} with $\S$ and $\T$ estimated using the method described in Proposition \ref{prop-estimating}; and OptShrink  \cite{nadakuditi2014optshrink}. The experiment is repeated $50$ times for each value of $\kappa$.

Figure \ref{fig-heter} shows the $\log_2$ mean relative errors of each method as a function of $\log_2(\kappa)$. For this model of $\S$ and $\T$, the condition number $\kappa$ is an increasing function of the parameter $\tau$ from Section \ref{sec-snr}. Consequently, Proposition \ref{prop-snr} suggests that whitening will improve the matrix SNR, and that the improvement should increase as $\kappa$ grows. This is precisely what Figure \ref{fig-heter} demonstrates; optimal spectral denoising with whitening by $\S$ and $\T$ does indeed outperform OptShrink, and the performance gap grows with $\kappa$. Using the estimated covariances, the performance is degraded but still outperforms OptShrink when $\kappa$ is large.

\subsection{Missing data}
\label{sec-numerics-missing}

We test  spectral denoising for  missing data (Algorithm \ref{alg-missing}) by comparing it to nuclear-norm regularized least-squares \cite{candes2010noise}, which estimates $\X$ by:
\begin{align}
\label{eq-xnuc}
\what \X^{\nuc} = \argmin_{\what \X \in \R^{p \times n}} 
    \left\{ \frac{1}{2}\|\F(\what \X) -  \y\|^2 
        + \theta \| \P^{1/2} \what \X \Q^{1/2} \|_* \right\}.
\end{align}
Here, $\|\cdot\|_*$ denotes the nuclear norm; $\F : \R^{p \times n} \to \R^M$ is the projection operator onto the $M$ observed samples; and $\P$ and $\Q$ are the diagonal matrices of sampling probabilities for rows and columns, respectively. We weight the nuclear norm by the square root of the sampling probabilities, as suggested in \cite{chen2015completing}. Following \cite{candes2010noise}, we choose the parameter $\theta$ so that when $\y$ is pure noise, $\what \X^{\nuc}$ is set to zero. It follows from the KKT conditions \cite{boyd2004convex} that this is equivalent to
\begin{math}
\theta = \|\P^{-1/2} \F^*(\y) \Q^{-1/2} \|_*,
\end{math}
which is approximated by $1 + \sqrt{\gamma}$. We solve \eqref{eq-xnuc} using the algorithm in \cite{ji2009agm}.

We generate a rank $r=5$ signal matrix $\X$ of size $p$-by-$n$, $p=200$, $n=400$, with singular values $\sqrt{\sqrt{\gamma} + 200k}$, $k = 1,\dots,5$, $\gamma=1/2$. Both the left and right singular vectors of $\X$ are uniformly random. We add to $\X$ a Gaussian noise matrix $\G$, where each entry has variance $\sigma^2/n$ for a specified value of $\sigma$. $\X + \G$ is then subsampled with row and column sampling probabilities each equispaced between $0.3$ and $0.7$. For each value of $\sigma$, the experiment is repeated 50 times. Figure \ref{fig-missing} displays the $\log_2$ mean relative errors. When $\sigma$ is large, spectral denoising is superior, whereas in the small $\sigma$ regime nuclear-norm regularized least-squares is better.
%%%This is consistent with the experimental results from \cite{dobriban2017optimal}.

%
%    table from /experiments/nongauss_table/table77.tex
%
\begin{table} 
\centering 
\begin{tabular}{|l| c  | c | c | c|}  
\hline  
 &  \multicolumn{4}{c|}{\textbf{Mean relative error}, $\mathbf{C}^\omega$} \\
\hline  
 $n$ & \textbf{Gaussian} & \textbf{Rademacher} & \textbf{t, df=10} & \textbf{t, df=3} \\
\hline  
$500$ & 2.956e-02 & 3.057e-02 & 2.928e-02 & 3.284e-01  \\
\hline  
$1000$ & 2.050e-02 & 2.127e-02 & 2.042e-02 & 4.208e-01  \\
\hline  
$2000$ & 1.459e-02 & 1.514e-02 & 1.466e-02 & 5.280e-01  \\
\hline  
$4000$ & 1.033e-02 & 1.019e-02 & 1.030e-02 & 6.535e-01  \\
\hline  
$8000$ & 7.459e-03 & 7.693e-03 & 7.495e-03 & 7.640e-01  \\
\hline  
\end{tabular}  

\vspace{\baselineskip}

%
%    table from /experiments/nongauss_table/table77.tex
%
\begin{tabular}{|l| c  | c | c | c|}  
\hline  
 &  \multicolumn{4}{c|}{\textbf{Max relative error}, $\mathbf{C}^\omega$} \\
\hline  
 $n$ & \textbf{Gaussian} & \textbf{Rademacher} & \textbf{t, df=10} & \textbf{t, df=3} \\
\hline  
$500$ & 1.117e-01 & 1.126e-01 & 1.241e-01 & 9.702e-01  \\
\hline  
$1000$ & 7.895e-02 & 8.411e-02 & 7.878e-02 & 9.837e-01  \\
\hline  
$2000$ & 6.520e-02 & 5.106e-02 & 6.907e-02 & 9.953e-01  \\
\hline  
$4000$ & 4.490e-02 & 4.056e-02 & 4.638e-02 & 9.950e-01  \\
\hline  
$8000$ & 3.035e-02 & 3.172e-02 & 3.174e-02 & 9.951e-01  \\
\hline  
\end{tabular}  
\caption{Average and maximum relative errors $\|\what \C^\omega - \C^\omega\|_{\Fr} / \|\C^\omega\|_{\Fr}$; see Section \ref{sec-nongaussian} for simulation details. For Gaussian, Rademacher, and t$_{10}$  distributions, the average errors decay approximately like $O(n^{-1/2})$.
%%%, as expected (see Remark \ref{rmk-accuracy}).
The errors for the t$_3$ distribution do not decay, indicating poor model fit.
}
\label{table-comega}
\end{table}

%
%    table from /experiments/nongauss_table/table77.tex
%
\begin{table} 
\centering 
\begin{tabular}{|l| c  | c | c | c|}  
\hline  
 &  \multicolumn{4}{c|}{\textbf{Mean relative error}, $\mathbf{D}$} \\
\hline  
 $n$ & \textbf{Gaussian} & \textbf{Rademacher} & \textbf{t, df=10} & \textbf{t, df=3} \\
\hline  
$500$ & 1.991e-02 & 2.039e-02 & 1.975e-02 & 1.757e-01  \\
\hline  
$1000$ & 1.414e-02 & 1.437e-02 & 1.412e-02 & 2.219e-01  \\
\hline  
$2000$ & 9.917e-03 & 1.015e-02 & 1.001e-02 & 2.788e-01  \\
\hline  
$4000$ & 7.027e-03 & 6.875e-03 & 7.005e-03 & 3.455e-01  \\
\hline  
$8000$ & 5.040e-03 & 5.224e-03 & 5.087e-03 & 4.141e-01  \\
\hline  
\end{tabular}  

\vspace{\baselineskip}

%
%    table from /experiments/nongauss_table/table77.tex
%
\begin{tabular}{|l| c  | c | c | c|}  
\hline  
 &  \multicolumn{4}{c|}{\textbf{Max relative error, $\mathbf{D}$}} \\
\hline  
 $n$ & \textbf{Gaussian} & \textbf{Rademacher} & \textbf{t, df=10} & \textbf{t, df=3} \\
\hline  
$500$ & 6.954e-02 & 7.491e-02 & 7.354e-02 & 9.663e-01  \\
\hline  
$1000$ & 5.044e-02 & 4.533e-02 & 5.505e-02 & 9.238e-01  \\
\hline  
$2000$ & 3.772e-02 & 3.578e-02 & 4.101e-02 & 9.529e-01  \\
\hline  
$4000$ & 2.468e-02 & 2.564e-02 & 2.712e-02 & 9.741e-01  \\
\hline  
$8000$ & 1.584e-02 & 1.836e-02 & 1.668e-02 & 9.779e-01  \\
\hline  
\end{tabular}  
\caption{Average and maximum relative errors $\|\what \D - \D\|_{\Fr} / \|\D\|_{\Fr}$; see Section \ref{sec-nongaussian} for simulation details. For Gaussian, Rademacher, and t$_{10}$  distributions, the average errors decay approximately like $O(n^{-1/2})$.
%%% as expected (see Remark \ref{rmk-accuracy}).
The errors for the  t$_3$ distribution do not decay, indicating poor model fit.
}
\label{table-d}
\end{table}

%

%
%   table from /experiments/rank_table/table77.tex
%
\begin{table} 
\centering 
\begin{tabular}{|l|| c  | c | c |}  
\hline  
 &  \multicolumn{3}{c|}{\textbf{Mean relative error}}  \\
\hline  
 \textbf{Noise type} & \textbf{Oracle} & \textbf{K-N} & \textbf{Naive}  \\
\hline  
Gaussian& 4.010e-01 & 4.010e-01 & 4.012e-01  \\
\hline  
Rademacher& 4.012e-01 & 4.012e-01 & 4.013e-01  \\
\hline  
t, df=10& 4.022e-01 & 4.022e-01 & 4.024e-01  \\
\hline  
t, df=5& 4.047e-01 & 4.086e-01 & 4.099e-01  \\
\hline  
t, df=4& 4.337e-01 & 4.484e-01 & 4.539e-01  \\
\hline  
t, df=3& 6.606e-01 & 7.459e-01 & 7.611e-01  \\
\hline  
t, df=2.5& 1.026e+00 & 1.287e+00 & 1.300e+00 \\
\hline  
\end{tabular}

\vspace{\baselineskip}

\begin{tabular}{|l|| c  | c | c |}  
\hline  
 &   \multicolumn{3}{c|}{\textbf{Max relative error}} \\
\hline  
 \textbf{Noise type}  & \textbf{Oracle} & \textbf{K-N} & \textbf{Naive} \\
\hline  
Gaussian & 4.420e-01 & 4.420e-01 & 4.421e-01 \\
\hline  
Rademacher & 4.459e-01 & 4.459e-01 & 4.459e-01 \\
\hline  
t, df=10 & 4.349e-01 & 4.349e-01 & 4.349e-01 \\
\hline  
t, df=5 & 1.280e+00 & 1.262e+00 & 1.262e+00 \\
\hline  
t, df=4 & 2.149e+00 & 2.113e+00 & 2.113e+00 \\
\hline  
t, df=3 & 4.922e+00 & 4.940e+00 & 4.943e+00 \\
\hline  
t, df=2.5 & 6.693e+00 & 6.703e+00 & 6.706e+00 \\
\hline  
\end{tabular}

\caption{Average and maximum relative errors of estimation; see Section \ref{sec-rank} for simulation details. The naive rank estimate $\hat r^{\naive}$ from \eqref{rank-naive} tends to overestimate the true rank $r=2$, whereas the estimate $\hat r^{\KN}$ of Kritchman and Nadler \cite{kritchman2008determining} is more accurate. However, the difference between the errors in the resulting estimates of $\X$ is not large. Both methods perform poorly for heavy tailed distributions.
}
\label{table-ranks1}
\end{table}

\begin{table} 
\centering 

\vspace{\baselineskip}
%
%   table from /experiments/rank_table/table77.tex
%
\begin{tabular}{|l|| c  | c | c |}  
\hline  
 &  \multicolumn{3}{c|}{\textbf{Mean rank}}  \\
\hline  
 \textbf{Noise type} & \textbf{Oracle} & \textbf{K-N} & \textbf{Naive}  \\
\hline  
Gaussian& 2.000e+00 & 2.000e+00 & 2.084e+00 \\
\hline  
Rademacher& 2.000e+00 & 2.000e+00 & 2.037e+00 \\
\hline  
t, df=10& 2.000e+00 & 2.000e+00 & 2.094e+00 \\
\hline  
t, df=5& 2.000e+00 & 2.128e+00 & 2.443e+00  \\
\hline  
t, df=4& 2.000e+00 & 2.858e+00 & 3.577e+00  \\
\hline  
t, df=3& 2.000e+00 & 7.875e+00 & 8.967e+00  \\
\hline  
t, df=2.5& 2.000e+00 & 1.643e+01 & 1.716e+01  \\
\hline  
\end{tabular}

\vspace{\baselineskip}
%
%   table from /experiments/rank_table/table77.tex
%
\begin{tabular}{|l|| c  | c | c |}  
\hline  
 &  \multicolumn{3}{c|}{\textbf{Max rank}} \\
\hline  
 \textbf{Noise type}  & \textbf{Oracle} & \textbf{K-N} & \textbf{Naive} \\
\hline  
Gaussian &       2&       2&       3 \\
\hline  
Rademacher &       2&       2&       3 \\
\hline  
t, df=10 &       2&       2&       3 \\
\hline  
t, df=5 &       2&       4&       5 \\
\hline  
t, df=4 &       2&       7&       8 \\
\hline  
t, df=3 &       2&      16&      17 \\
\hline  
t, df=2.5 &       2&      26&      26 \\
\hline  
\end{tabular}

\caption{Average and maximum rank estimates; see Section \ref{sec-rank} for simulation details. The naive rank estimate $\hat r^{\naive}$ from \eqref{rank-naive} tends to overestimate the true rank $r=2$, whereas the estimate $\hat r^{\KN}$ of Kritchman and Nadler \cite{kritchman2008determining} is more accurate. However, the difference between the errors in the resulting estimates of $\X$ is not large. Both methods perform poorly for heavy tailed distributions.
}
\label{table-ranks2}
\end{table}

\subsection{Non-Gaussian noise}
\label{sec-nongaussian}

The optimal spectral denoiser requires estimation of the weighted inner product matrices $\D$, $\wtilde \D$, $\C^\omega$ and $\wtilde \C^\omega$. The formulas for the entries of these matrices provided by Theorem \ref{thm-asymptotics} assumes that the noise matrix $\G$ is Gaussian. However, it is of interest whether the same formulas may be applied to non-Gaussian noise. To partially address this question, we compare the finite sample accuracy of the formulas in Theorem \ref{thm-asymptotics} for different noise distributions.

For each  $n$, we generate $\Y = \X + \G$ of size $p$-by-$n$, where $p = 2n$. The signal has rank $r=2$, with singular values $\gamma^{1/4} + 2$ and $\gamma^{1/4}+ 3$; $\u_1$ is uniformly equal to $1/\sqrt{p}$, and $\u_2$ is $1/\sqrt{p}$ on entries $1,\dots,p/2$, and $-1/\sqrt{p}$ on entries $p/2+1,\dots,p$. $\v_1$ and $\v_2$ are generated similarly, with $n$ in place of $p$. The noise matrix has iid entries of variance $1/n$, drawn from a specified distribution: Gaussian, Rademacher, t$_{10}$ or t$_3$, where the t distributions are normalized to have variance $1/n$.

%%%The t$_3$ distribution may be considered ``heavy tailed'' as it has fewer than three finite moments.

The $p$-by-$p$ weight matrix $\Omega$ is diagonal with diagonal entries $1,\dots,3p/4$ equal to $1$, and the remaining entries $0$. We evaluate the true matrix $\E$ and use formulas \eqref{fmla-comega} and \eqref{fmla-d} to predict $\D$ and $\C^{\omega}$. For each draw, we compute the actual inner product matrices $\what \D$ and $\what \C^\omega$ using the left singular vectors $\hat \u_1$ and $\hat \u_2$ of $\Y$. Due to the ambiguity in signs, we make all entries of the matrices positive. We then compute the relative errors $\|\what \D - \D\|_{\Fr} / \|\D\|_{\Fr}$ and $\|\what \C^\omega - \C^\omega\|_{\Fr} / \|\C^\omega\|_{\Fr}$.

For each noise type and each value of $n = 500k$, $k=1,2,4,8,16$, the experiment is repeated 1000 times. The average and maximum relative errors are recorded in Table \ref{table-comega} for $\C^{\omega}$, and in  Table \ref{table-d} for $\D$. Both the average and maximum  errors for Gaussian noise very nearly match those for the Rademacher and t$_{10}$ distributions. However, the errors for the heavier tailed t$_3$ distribution are much larger, indicating that the theory breaks down for this noise model. The errors for the Gaussian, Rademacher, and t$_{10}$ distributions appear to decay approximately like $O(n^{-1/2})$; this is the rate we expect from \cite{bao2018singular} and Theorem 2.19 in \cite{benaych2012singular}. The errors for the t$_3$ distribution do not exhibit such decay with increasing $n$, indicating that the model does not match.

\subsection{Rank estimation}
\label{sec-rank}

In this section, we explore estimation of the rank $r$ of $\X$ from the observed matrix $\Y$, a topic that has received considerable attention \cite{kritchman2008determining, kritchman2009non, dobriban2017deterministic, dobriban2017factor, passemier2014estimation, passemier2012determining}. The naive estimator $\hat r^{\naive}$ is defined by
\begin{align}
\label{rank-naive}
\hat r^{\naive} = \#\{k : \lambda_k > 1 + \sqrt{\gamma}\};
\end{align}
this simply counts the number of $\Y$'s singular values exceeding $1 + \sqrt{\gamma}$, the asymptotically largest singular value of the noise matrix $\G$. It is known that $\hat r^{\naive}$ may overestimate the rank; see, e.g., \cite{johnstone2018tail}. The rank estimator of Kritchman and Nadler from \cite{kritchman2008determining}, denoted by $\hat r^{\KN}$, is designed to prevent attributing noisy singular values to signal. We compare the performance of $\hat r^{\naive}$ and $\hat r^{\KN}$ for different noise distributions in terms of the accuracy of estimating $r$ and the effect on the denoising error.

For $p=300$ and $n=600$, we generate a $p$-by-$n$ signal matrix $\X$ with rank $r=2$ and singular values $\gamma^{1/4} + 1$ and $\gamma^{1/4}+ 2$. The left singular vector $\u_1$ is uniformly equal to $1/\sqrt{p}$, and $\u_2$ is $1/\sqrt{p}$ on entries $1,\dots,p/2$, and $-1/\sqrt{p}$ on entries $p/2+1,\dots,p$. The right singular vectors $\v_1$ and $\v_2$ are generated the same way, with $n$ in place of $p$. The noise matrix $\G$ had iid entries of variance $1/n$, drawn from one a specified distribution, namely: Gaussian, Rademacher, or the t distributions with $10$, $5$, $4$, $3$ and $2.5$ degrees of freedom, where the t distributions are normalized to have variance $1/n$. The $p$-by-$p$ weight matrix $\Omega$ is diagonal with diagonal entries linearly spaced between $1/p$ and $1$. The $n$-by-$n$ weight matrix $\Pi$ is also diagonal, with diagonal entries linearly spaced between $1/p$ and $1/\gamma$. In each run, we apply Algorithm \ref{alg:denoise} with the oracle $r=2$, the naive $\hat r^{\naive}$ from \eqref{rank-naive}, and $\hat r^{\KN}$ from \cite{kritchman2008determining}, with $0.1$ confidence level\footnote{The code for computing $\hat r^{\KN}$ was taken from Boaz Nadler's website: \url{www.wisdom.weizmann.ac.il/~nadler/Rank_Estimation/rank_estimation.html}}. For each noise distribution, the experiment is repeated $1000$ times. In Tables \ref{table-ranks1} and \ref{table-ranks2} we record the relative errors $\|\Omega(\what \X - \X)\Pi^T\|_{\Fr} / \|\Omega \X\Pi^T\|_{\Fr}$ and the estimated ranks.

For the Gaussian, Rademacher, and t$_{10}$ distributions, the Kritchman-Nadler estimate $\hat r^{\KN}$ is typically closer to the true rank, $r=2$, than is the naive estimate $\hat r^{\naive}$. However, the average and maximum errors are close regardless of the rank estimator used, since even when $\hat r^{\naive}=3$, the third singular value of $\Y$ is so close to the detection edge $1 + \sqrt{\gamma}$ that the estimates of the cosines $c_3$ and $\tilde c_3$ are nearly $0$. With the t$_5$ distribution, both $\hat r^{\naive}$ and $\hat r^{\KN}$ are more likely to overestimate the true rank.  While the average errors are close to those for the Gaussian, Rademacher, and t$_{10}$ distributions, the maximum errors are much larger, indicating that while this noise distribution's ``typical'' behavior may be close to the thinner tailed ones, a small number of extreme draws of $\G$ can result in very poor performance. For the t distributions with $4$, $3$, and $2.5$ degrees of freedom, both $\hat r^{\naive}$ and $\hat r^{\KN}$ drastically overestimate the rank, and the resulting relative errors are enormous.

\section{Conclusion}
\label{sec-conclusion}

This paper has introduced a family of spectral denoisers for low-rank matrix estimation, which generalizes singular value shrinkage. We have derived optimal spectral denoisers for weighted loss functions, and discussed applications. By judiciously combining these denoisers for different weights we contructed the method of localized denoising, which outperforms singular value shrinkage under heterogeneity. While this paper has focused on theoretical and algorithmic development, in future work we plan to apply the methods to problems where related but suboptimal methods have previously been employed. This includes the problems of denoising and deconvolution of images from cryoelectron microscopy \cite{bhamre2016denoising}; 3-D reconstruction of heterogeneous molecules from noisy images \cite{anden2018structural}; and denoising XFEL images \cite{liu2016epca, zhao2018steerable}.

%%% We also showed how spectral denoising with weighted loss may be applied to submatrix denoising, denoising matrices with doubly-heteroscedastic noise, and denoising with missing data.

%%%The theory we have developed holds for Gaussian noise. However, we demonstrated numerically that the predicted asymptotic behavior appears to hold for other thin tailed noise distributions. It is of considerable interest to extend the theoretical results to this more general setting. This is particularly true for problems with missing data, as well as other instances of the linearly-transformed spiked model \cite{dobriban2017optimal}.

\section*{Acknowledgements}

I thank Elad Romanov and Amit Singer for stimulating discussions related to this work, Edgar Dobriban for valuable feedback on an earlier version of this manuscript, and the reviewers for their helpful comments. I acknowledge support from the NSF BIGDATA award IIS 1837992 and BSF award 2018230.

\bibliographystyle{plain}
\bibliography{refs_weighted}

\appendix

\section{Proof of Theorem \ref{thm-asymptotics}}
\label{appendix-asymptotics}

The proof of Theorem \ref{thm-asymptotics} is similar to the analysis found in \cite{leeb2019optimal}, in that it rests on the same decomposition of the empirical singular vectors $\hat \u_j$ and $\hat \v_j$ into the signal and residual components. If $\a$ and $\b$ are vectors of the same dimension, we will write $\a \sim \b$ as a short-hand for $\|\a - \b\| \to 0$ almost surely as $p,n \to \infty$. The statements are symmetric in the left and right singular vectors, so for compactness we will only prove them for the left ones. The proofs for the other side are identical.

Because the noise matrix $\G$ has an isotropic distribution, we can write:
\begin{align}
\label{eq-resids}
\hat \u_j \sim c_j \u_j + s_j \tilde \u_j,
\end{align}
where $\tilde \u_j$ is a unit vector that is uniformly random over the sphere in the subspace orthogonal to $\u_1,\dots,\u_r$ (see \cite{paul2007asymptotics}). Because $\tilde \u_j$ is uniformly random, it is asymptotically orthogonal to any independent unit vector $\w$; that is,
\begin{align}
\label{eq-orth}
\tilde \u_j^T \w \sim 0.
\end{align}
Furthermore, $\tilde \u_j$ satisfies the normalized trace formula, namely if $\A$ is any matrix with bounded operator norm, then
\begin{align}
\label{eq-trace}
\tilde \u_j^T \A \tilde \u_j \sim \frac{1}{p}\tr{\A}.
\end{align}
We refer the reader to \cite{benaych2011fluctuations, hanson1971bound, wright1973bound, rudelson2013hanson} for details. We will use \eqref{eq-orth} and \eqref{eq-trace} repeatedly. Furthermore, when $j \ne k$ it follows from Lemma A.2 in \cite{leeb2019optimal} that
\begin{align}
\label{eq-orth2}
\tilde \u_j^T \A \tilde \u_k \sim 0.
\end{align}

Applying $\Omega$ to each side of \eqref{eq-resids}, we have:
\begin{align}
\label{eq-main}
\Omega \hat \u_j \sim c_j \Omega \u_j + s_j \Omega \tilde \u_j.
\end{align}
The proofs of the identities in Theorem \ref{thm-asymptotics} now follow by manipulating the asymptotic equation \eqref{eq-main} appropriately, in conjunction with \eqref{eq-orth}, \eqref{eq-trace} and \eqref{eq-orth2}.

We first show the formulas for $c_{jk}^\omega$. We take inner products of each side of \eqref{eq-main} with $\Omega \u_k$:
\begin{align}
c_{jk}^\omega \sim \la \Omega \hat \u_j , \Omega \u_k \ra
\sim c_j \la  \Omega \u_j, \Omega \u_k \ra 
    + s_j \la \Omega \tilde \u_j,  \Omega \u_k \ra
\sim c_j \la  \Omega \u_j, \Omega \u_k \ra
\sim c_j e_{jk},
\end{align}
where we have used \eqref{eq-orth}.

To derive the formula for $d_{j}$, we take the squared norm of each side of \eqref{eq-main}:
\begin{align}
d_j \sim
\|\Omega \hat \u_j\|^2 
    \sim c_j^2\|\Omega \u_j\|^2  + s_j^2 \|\Omega \tilde \u_j\|^2
    \sim c_j^2\alpha_j + s_j^2 \mu.
\end{align}
The first asymptotic equivalence follows from \eqref{eq-orth}, and the second from \eqref{eq-trace}.

Finally, we derive the formula for $d_{jk}$, $j \ne k$. From \eqref{eq-main}, we have
\begin{align}
\langle \Omega \hat \u_j , \Omega \hat \u_k \rangle
&\sim c_j c_k \langle \Omega \u_j , \Omega \u_k \rangle
    + s_j s_k \la \Omega \tilde \u_j, \Omega \tilde \u_k \ra
    + s_j c_k \la \Omega \tilde \u_j, \Omega \u_k\ra
    + c_j s_k \la \Omega \u_j, \Omega \tilde \u_k\ra.
\end{align}
From \eqref{eq-orth} and \eqref{eq-orth2}, the terms involving $\tilde \u_j$ and $\tilde \u_k$ vanish, and we are left with
\begin{align}
d_{jk} \sim
\langle \Omega \hat \u_j , \Omega \hat \u_k \rangle
\sim c_j c_k \langle \Omega \u_j , \Omega \u_k \rangle
\sim c_j c_k e_{jk} .
\end{align}
This completes the proof of Theorem \ref{thm-asymptotics}.

\section{Proof of Theorem \ref{thm-ls}}
\label{appendix-ls}

The target matrix $\X$ may be written
\begin{align}
\X = \sum_{k=1}^{r} t_k \u_k \v_k^T
= \U \diag(\t) \V^T,
\end{align}
and our estimate $\what \X$ is of the form
\begin{align}
\what \X = \what \U \what \B \what \V^T,
\end{align}
where $\U = [\u_1,\dots,\u_r]$, $\V = [\v_1,\dots,\v_r]$, $\what \U = [\hat \u_1,\dots,\hat \u_r]$, $\what \V = [\hat \v_1,\dots,\hat\v_r]$, and $\t = (t_1,\dots,t_r)^T$.

Define $\W = \Omega \U$, $\Z = \Pi \V$, $\what \W = \Omega \what \U$, and $\what \Z = \Pi \what \V$. We may then write the weighted loss as follows:
\begin{align}
\L(\what \X, \X)
= \| \Omega (\what \X - \X)\Pi^T \|_{\Fr}^2
= \| \Omega \what \X \Pi^T - \Omega\X\Pi^T \|_{\Fr}^2
= \| \what \W \what \B \what \Z^T - \W \diag(\t) \Z^T \|_{\Fr}^2,
\end{align}
which is the \emph{unweighted} Frobenius loss between $\what \W \what \B \what \Z^T$ and $\W \diag(\t) \Z^T$. Continuing, we have:
\begin{align}
\label{eq3400}
\L(\what \X, \X)
&= \| \W \diag(\t) \Z^T - \what \W \what \B \what \Z^T\|_{\Fr}^2
\nonumber \\
&=  \| \W \diag(\t) \Z^T \|_{\Fr}^2 + \| \what \W \what \B \what \Z^T\|_{\Fr}^2
    - 2 \la \W \diag(\t) \Z^T , \what \W \what \B \what \Z^T \ra_{\Fr}
\nonumber \\
&=  \la  \W^T \W \diag(\t) \Z^T \Z , \diag(\t)\ra_{\Fr}
    + \la  \what \W^T \what \W \what \B \what \Z^T \what \Z , \what \B\ra_{\Fr}
    - 2 \la \what \W^T \W \diag(\t) \Z^T  \what \Z,  \what \B \ra_{\Fr}
\nonumber \\
&\sim  \la   \E \diag(\t) \wtilde \E, \diag(\t) \ra_{\Fr}
    + \la  \D \what \B \wtilde \D , \what \B\ra_{\Fr}
    - 2 \la \C \diag(\t) \wtilde \C^T,  \what \B \ra_{\Fr}.
\end{align}
%
%%%where $\E = \W^T \W$, $\wtilde \E = \Z^T \Z$, $\D = \what \W^T \what \W$, $\wtilde \D = \what \Z^T \what \Z$, $\C = \what \W^T \W$, and $\wtilde \C = \what \Z^T \Z$. Written in this way, we see that the objective $\|\A - \what \A\|_{\Fr}^2$ depends only on the inner product matrices $\E$, $\wtilde \E$, $\D$, $\wtilde \D$, $\C$, and $\wtilde \C$.

Defining the operator $\TT$ by $\TT(\what \B) = \D \what \B \wtilde \D$, the pseudoinverse of $\TT$ is given by $\TT^+(\B) = \D^+ \B \wtilde \D^+$. Consequently, the choice of $\what \B$ that minimizes $\L(\what \X , \X)$ is given by:
\begin{align}
\what \B = \D^+ \C \diag(\t) \wtilde \C^T \wtilde \D^+.
\end{align}
The error may then be evaluated by substituting this expression for $\what \B$ into \eqref{eq3400}, completing the proof.

\section{Proof of Theorem \ref{thm-diag}}
\label{appendix-diag}

Under weighted orthogonality, $e_{jk} = \tilde e_{jk} = 0$ whenever $j \ne k$, and so $d_{jk} = \tilde d_{jk} = c_{jk}^\omega  = \tilde c_{jk}^\omega = 0$ when $j \ne k$ as well. Consequently, the matrices $\E$, $\wtilde \E$, $\D$, $\wtilde \D$, $\C$, and $\wtilde \C$ are diagonal. The optimal $\what \B$ is given by:
\begin{align}
\what \B = \D^+ \C \diag(\t) \wtilde \C \wtilde \D^+,
\end{align}
which is also diagonal, with diagonal entries
\begin{align}
\hat t_k 
= \frac{t_k c_k^\omega \tilde c_k^\omega}{d_k \tilde d_k}
= \frac{t_k c_k \alpha_k  \tilde c_k \beta_k}
    {(c_k^2 \alpha_k + s_k^2 \mu)(\tilde c_k^2 \beta_k + \tilde s_k^2 \nu)}
= t_k c_k  \tilde c_k\frac{\alpha_k \beta_k}
    {(c_k^2 \alpha_k + s_k^2 \mu)(\tilde c_k^2 \beta_k + \tilde s_k^2 \nu)},
\end{align}
which is the desired expression.

\section{Proof of Proposition \ref{prop-shrink}}
\label{appendix-shrink}

Suppose a coordinate has signal strength $t = t_k$ (we drop the subscript as we are only considering one component). We may assume without loss of generality (and by rescaling $\alpha$ and $\beta$) that $\mu = \nu = 1$. Consequently, the optimal singular value is equal to:
\begin{align}
\hat t = t c \tilde c \cdot \frac{\alpha}{c^2 \alpha  + s^2} 
    \cdot \frac{\beta}{\tilde c^2 \beta + \tilde s^2}.
\end{align}
By taking $\alpha$ and $\beta$ sufficiently large, this value can be made arbitrarily close to
\begin{align}
\frac{t}{c \tilde c}
= \frac{t \sqrt{(1+\gamma/t^2)(1+1/t^2)} }{1 - \gamma/t^4}
= \frac{\lambda}{1 - \gamma/t^4}
> \lambda.
\end{align}
That is, the optimal singular value $\hat t$ will be greater than the observed singular value $\lambda$ in this parameter regime.

On the other hand, if $\beta \le 1 = \nu$, we have:
\begin{align}
\frac{\hat t}{\lambda} = \frac{1}{\lambda} t c \tilde c \cdot \frac{\alpha}{\alpha c^2 + s^2} 
    \cdot \frac{\beta}{\beta \tilde c^2 + \tilde s^2}.
\le \frac{1}{\lambda} t \frac{\tilde c}{c}
= \frac{t}{\sqrt{(t^2+1)(t^2+\gamma)}} t \sqrt{ \frac{t^2 + \gamma}{t^2 + 1} }
= \frac{t^2}{t^2+1}
\le 1,
\end{align}
which shows that $\hat t \le \lambda$. A nearly identical proof works if $\alpha \le \mu$. This completes the proof.

\section{Proof of Proposition \ref{prop-monotone}}
\label{appendix-monotone}

Without loss of generality, we will assume $\mu = \nu = 1$. We consider the functions $c(t) = \sqrt{(1 - \gamma/t^4) / (1 + \gamma/t^2)}$ and $\tilde c(t) = \sqrt{(1 - \gamma/t^4) / (1 + 1/t^2)}$. Define the functions $\vphi(t)$ and $\psi(t)$ by
\begin{align}
\vphi(t) = \frac{\alpha c(t)}{\alpha c(t)^2 + 1 - c(t)^2}
\end{align}
and
\begin{align}
\psi(t) = \frac{\beta \tilde c(t)}{\beta \tilde c(t)^2 + 1 - \tilde c(t)^2}.
\end{align}
Then we may write the optimal singular value $\hat t$ as a function $f(t)$ as follows:
\begin{align}
f(t) = t \vphi(t) \psi(t).
\end{align}
Let us assume that $\alpha \le 1$; the proof for $\beta \le 1$ will be nearly identical. We wish to show that $f'(t) \ge 0$, for $t > \gamma^{1/4}$. We have
\begin{align}
\frac{f'(t)}{f(t)} = \frac{\vphi'(t)}{\vphi(t)} + \frac{\psi'(t)}{\psi(t)} + \frac{1}{t},
\end{align}
and since $f(t) > 0$, we must show that the right side is positive. It is straightforward to verify that
\begin{align}
\vphi'(t)
= \frac{\alpha c'(t) [1 - (\alpha - 1)c(t)^2]}{[1 + (\alpha - 1)c(t)^2]^2}
\end{align}
from which it follows that
\begin{align}
\frac{\vphi'(t)}{\vphi(t)}
= \frac{c'(t)}{c(t)} \frac{1 - (\alpha-1)c(t)^2}{1 + (\alpha-1)c(t)^2}
\ge \frac{c'(t)}{c(t)}.
\end{align}
Similarly, we can show
\begin{align}
\frac{\psi'(t)}{\psi(t)}
= \frac{\tilde c'(t)}{\tilde c(t)} 
    \frac{1 - (\beta-1)\tilde c(t)^2}{1 + (\beta-1)\tilde c(t)^2}
\ge - \frac{\tilde c'(t)}{\tilde c(t)} .
\end{align}
Consequently, it is enough to show
\begin{align}
\label{ineq73200}
\frac{c'(t)}{c(t)} - \frac{\tilde c'(t)}{\tilde c(t)} + \frac{1}{t} \ge 0.
\end{align}
Direction computation shows
\begin{align}
\label{c-ratio}
\frac{c'(t)}{c(t)} = \gamma \frac{t^4 + 2 t^2 + \gamma}{t (t^2+\gamma)(t^4-\gamma)}
\end{align}
and
\begin{align}
\label{ctilde-ratio}
\frac{\tilde c'(t)}{\tilde c(t)} 
= \frac{t^4 + 2 \gamma t^2 + \gamma}{t (t^2+1)(t^4-\gamma)}.
\end{align}
Substituting \eqref{c-ratio} and \eqref{ctilde-ratio} into the left side of \eqref{ineq73200} and multiplying by $t (t^2 + \gamma) (t^2 + 1) $, we get:
\begin{align}
& t (t^2 + \gamma) (t^2 + 1)
    \left(\frac{c'(t)}{c(t)} - \frac{\tilde c'(t)}{\tilde c(t)} + \frac{1}{t}\right)
\nonumber \\
& = t (t^2 + \gamma) (t^2 + 1)
    \left(\gamma \frac{t^4 + 2 t^2 + \gamma}{t (t^2+\gamma)(t^4-\gamma)} 
    - \frac{t^4 + 2 \gamma t^2 + \gamma}{t (t^2+1)(t^4-\gamma)} 
    + \frac{1}{t}\right)
\nonumber \\
& = t^4 + 2\gamma t^2 + \gamma > 0,
\end{align}
which is the desired result.

\section{Proof of Theorem \ref{thm-compare1}}
\label{appendix-compare1}

We denote by $\hat t_1^{\shr},\dots,\hat t_r^{\shr}$ the singular values of $\what \X^{\shr}$, and $\hat \t^{\shr} = (\hat t_1^{\shr},\dots,\hat t_r^{\shr})^T$. We may then write
\begin{align}
\what \X^{\shr} = \what \U \diag(\hat \t^{\shr}) \what \V^T.
\end{align}
This is a spectral denoiser (in the set $\SS$), and hence its weighted loss with weights $\Omega_i$ and $\Pi_j$ cannot be less than that of the optimal spectral denoiser $\what \X_{(i,j)}^{\loc}$. That is,
\begin{align}
\| \Omega_i (\what \X_{(i,j)}^{\loc}- \X) \Pi_j^T \|_{\Fr}^2
\le \| \Omega_i (\what \X^{\shr} - \X) \Pi_j^T \|_{\Fr}^2
\end{align}

Because the $\Omega_i$ and $\Pi_j$ are pairwise orthogonal projections which sum to the identity, the total Frobenius loss can be decomposed:
\begin{align}
\| \what \X^{\loc}- \X \|_{\Fr}^2
&= \sum_{i=1}^{I} \sum_{j=1}^{J} \| \Omega_i (\what \X^{\loc}- \X) \Pi_j^T \|_{\Fr}^2
= \sum_{i=1}^{I} \sum_{j=1}^{J} \| \Omega_i (\what \X_{(i,j)}^{\loc}- \X) \Pi_j^T \|_{\Fr}^2
\nonumber \\
&\le \sum_{i=1}^{I} \sum_{j=1}^{J} \| \Omega_i (\what \X^{\shr} - \X) \Pi_j^T \|_{\Fr}^2
= \| \what \X^{\shr}- \X \|_{\Fr}^2,
\end{align}
which is the desired inequality.

\section{Proof of Theorem \ref{thm-compare2}}
\label{appendix-compare2}

For $1 \le k \le r$,  $1 \le i \le I$, and $1 \le j \le J$, let $\alpha_k^{(i)} = \| \Omega_i \u_k\|^2$, $\mu^{(i)} = \tr{\Omega_i} / p$, $\beta_k^{(j)} = \| \Pi_j \v_k\|^2$, and $\nu^{(j)} = \tr{\Pi_j} / n$. Then 
\begin{align}
\sum_{i=1}^{I} \alpha_k^{(i)}
= \sum_{i=1}^{I} \mu^{(i)}
= \sum_{j=1}^{J} \beta_k^{(j)}
= \sum_{j=1}^{J} \nu^{(j)}
= 1.
\end{align}

Let $\what \X_{(i,j)}^{\dd}$ be the optimal diagonal denoiser with weights $\Omega_i$ and $\Pi_j$. Because of the weighted orthogonality condition, Theorem \ref{thm-diag} states that the AMSE for $\what \X_{(i,j)}^{\dd}$ is
\begin{align}
\|\Omega_i(\what \X_{(i,j)}^{\dd} - \X )\Pi_j^T\|_{\Fr}^2
&= \sum_{k=1}^{r} t_k^2 \alpha_k^{(i)} \beta_k^{(j)} \left( 1 
    - c_k^2 \tilde c_k^2 \cdot \frac{\alpha_k^{(i)}}{c_k^2\alpha_k^{(i)} + s_k^2 \mu^{(i)}}
    \cdot \frac{ \beta_k^{(j)}}{\tilde c_k^2\beta_k^{(j)} + \tilde s_k^2 \nu^{(j)}}\right).
\end{align}
Since $\what \X_{(i,j)}^{\loc}$ minimizes the weighted error with weights $\Omega_i$ and $\Pi_j$, we have:
\begin{align}
\label{err-loc4700}
\|\what \X^{\loc} - \X \|_{\Fr}^2
&=\sum_{i=1}^{I} \sum_{j=1}^{J} \|\Omega_i(\what \X^{\loc} - \X)\Pi_j^T \|_{\Fr}^2
=\sum_{i=1}^{I} \sum_{j=1}^{J} \|\Omega_i(\what \X_{(i,j)}^{\loc} - \X)\Pi_j^T \|_{\Fr}^2
\nonumber \\
&\le\sum_{i=1}^{I} \sum_{j=1}^{J} \|\Omega_i(\what \X_{(i,j)}^{\dd} - \X)\Pi_j^T \|_{\Fr}^2
\nonumber \\
&=\sum_{i=1}^{I} \sum_{j=1}^{J} 
    \sum_{k=1}^{r} t_k^2 \alpha_k^{(i)} \beta_k^{(j)} \left( 1 
    - c_k^2 \tilde c_k^2 \cdot \frac{\alpha_k^{(i)}}{c_k^2\alpha_k^{(i)} + s_k^2 \mu^{(i)}}
    \cdot \frac{ \beta_k^{(j)}}{\tilde c_k^2\beta_k^{(j)} + \tilde s_k^2 \nu^{(j)}}\right).
\nonumber \\
&= \sum_{k=1}^{r} t_k^2 \left( 1  - 
     c_k^2 \tilde c_k^2 \sum_{i=1}^{I} \sum_{j=1}^{J} 
        \frac{(\alpha_k^{(i)})^2}{c_k^2\alpha_k^{(i)} + s_k^2 \mu^{(i)}}
        \cdot \frac{ (\beta_k^{(j)})^2 }{\tilde c_k^2\beta_k^{(j)} + \tilde s_k^2 \nu^{(j)}}
    \right).
\end{align}

On the other hand, the error obtained by $\what \X^{\shr}$ is equal to
\begin{align}
\label{err-shr4700}
\|\what \X^{\shr} - \X \|_{\Fr}^2 = \sum_{k=1}^{r} t_k^2 (1 - c_k^2 \tilde c_k^2).
\end{align}

Comparing \eqref{err-loc4700} and \eqref{err-shr4700}, the result will follow if we can show that for each $1 \le k \le r$,
\begin{align}
\sum_{i=1}^{I} \sum_{j=1}^{J} 
        \frac{(\alpha_k^{(i)})^2}{c_k^2\alpha_k^{(i)} + s_k^2 \mu^{(i)}}
        \cdot \frac{ (\beta_k^{(j)})^2 }{\tilde c_k^2\beta_k^{(j)} + \tilde s_k^2 \nu^{(j)}}
\ge 1,
\end{align}
where the inequality is strict so long as one of $\u_k$ or $\v_k$ is not generic with respect to some $\Omega_i$ or $\Pi_j$; or equivalently, either $\alpha_k^{(i)} \ne \mu^{(i)}$ for some $i$, or $\beta_k^{(j)} \ne \nu^{(j)}$ for some $j$. Because
\begin{align}
\sum_{i=1}^{I} \sum_{j=1}^{J} 
        \frac{(\alpha_k^{(i)})^2}{c_k^2\alpha_k^{(i)} + s_k^2 \mu^{(i)}}
        \cdot \frac{ (\beta_k^{(j)})^2 }{\tilde c_k^2\beta_k^{(j)} + \tilde s_k^2 \nu^{(j)}}
=
\left(\sum_{i=1}^{I} \frac{(\alpha_k^{(i)})^2}{c_k^2\alpha_k^{(i)} + s_k^2 \mu^{(i)}}
    \right) \cdot
\left(\sum_{j=1}^{J} 
    \cdot \frac{ (\beta_k^{(j)})^2 }{\tilde c_k^2\beta_k^{(j)} + \tilde s_k^2 \nu^{(j)}}
    \right),
\end{align}
it is enough to show that
\begin{align}
\sum_{i=1}^{I} \frac{(\alpha_k^{(i)})^2}{c_k^2\alpha_k^{(i)} + s_k^2 \mu^{(i)}}
\ge 1,
\end{align}
with the inequality being strict so long as $\alpha_k^{(i)} \ne \mu^{(i)}$ for some $i$.

For each $1 \le i \le I$, let $r_i = \alpha_k^{(i)} / \mu^{(i)}$. Then
\begin{align}
\sum_{i=1}^{I} \frac{(\alpha_k^{(i)})^2}{c_k^2\alpha_k^{(i)} + s_k^2 \mu^{(i)}}
= \sum_{i=1}^{I} \mu^{(i)} \frac{r_i^2}{c_k^2 r_i + s_k^2}.
\end{align}

The function $F(r) = r^2 / (c_k^2 r + s_k^2)$ is convex. Since $\sum_{i=1}^{I} \mu^{(i)} = 1$, Jensen's inequality implies
\begin{align}
\sum_{i=1}^{I} \mu^{(i)} \frac{r_i^2}{c_k^2 r_i + s_k^2}
= \sum_{i=1}^{I} \mu^{(i)} F(r_i)
\ge F\left( \sum_{i=1}^{I} \mu^{(i)} r_i\right)
= F\left( \sum_{i=1}^{I} \alpha_k^{(i)} \right)
= F(1) = 1,
\end{align}
which is the desired inequality. The inequality will be strict so long as $r_i = \alpha_k^{(i)}/ \mu^{(i)}$ is not constantly equal to $1$ over $i$, or equivalently if $\alpha_k^{(i)} \ne \mu^{(i)}$ for some $i$. This is the desired result.

\section{Proof of Proposition \ref{prop-merge}}
\label{appendix-merge}

Since $\Y_0 = \Omega \Y \Pi^T$ has only $n_0$ columns, to ensure that the scaling of the noise matches that of the standard spiked model, we must multiply it by $\sqrt{n / n_0} = 1/\sqrt{\nu}$. We define $\wtilde \Y_0 = \Y_0 / \sqrt{\nu}$ and $\wtilde \X_0 = \X_0 / \sqrt{\nu}$. Then $\wtilde \Y_0$ follows a standard spiked model with signal matrix $\wtilde \X_0$.

For $1 \le k \le r$, we let $\u_k$ and $\v_k$ denote the $k^{th}$ singular vectors of $\X$; $\hat \u_k$ and $\hat \v_k$ denote the $k^{th}$ singular vectors of $\Y$; $\u_k^0$ and $\v_k^0$ denote the $k^{th}$ singular vectors of $\X_0$ (and $\wtilde \X_0$); and $\hat \u_k^0$ and $\hat \v_k^0$ denote the $k^{th}$ singular vectors of $\Y_0$ (and $\wtilde \Y_0$). We let $t_1^0,\dots,t_r^0$ denote the singular values of $\wtilde \X_0$. We also let $\gamma_0 = p_0/n_0 = (\mu / \nu) \gamma$ be the aspect ratio of the submatrix.

If $t_1,\dots,t_r$ are the singular values of the full $p$-by-$n$ signal matrix $\X$, then we may write the rescaled submatrix $\wtilde \X_0$ as
\begin{align}
\label{X0-svd}
\wtilde \X_0 = \Omega\X\Pi^T / \sqrt{\nu} 
= \frac{1}{\sqrt{\nu}}\sum_{k=1}^{r} t_k \Omega \u_k \v_k^T \Pi^T
= \sum_{k=1}^{r} t_k \sqrt{\frac{\alpha_k \beta_k}{\nu}}
    \frac{\Omega \u_k}{\|\Omega \u_k\|} \left(\frac{\Pi \v_k}{\|\Pi \v_k\|} \right)^T.
\end{align}
Because the $\Omega \u_k$ and $\Pi \v_k$ are assumed to by pairwise orthogonal, \eqref{X0-svd} is the SVD of $\wtilde \X_0$. Consequently:
\begin{align}
\u_k^0 = \frac{\Omega \u_k}{\|\Omega \u_k\|},
\quad \v_k^0 = \frac{\Pi \v_k}{\|\Pi \v_k\|},
\quad t_k^0 = t_k \sqrt{\frac{\alpha_k \beta_k}{\nu}}.
\end{align}
We define the cosines
\begin{align}
c_k^0 = \la \hat \u_k^0, \u_k^0 \ra,
\quad \tilde c_k^0 = \la \hat \v_k^0, \v_k^0 \ra.
\end{align}
Following Remark \ref{rmk-positive}, we may assume that the singular vectors have been chosen so that both $c_k^0$ and $\tilde c_k^0$ are non-negative. Then the AMSE obtained by first applying optimal singular value shrinkage to $\wtilde \Y_0$, and then rescaling by $\nu$, is
\begin{align}
\label{amse-Xshr}
\| \what \X_0^{\shr} - \X_0\|_{\Fr}^2 
\sim \nu \sum_{k=1}^{r} (t_k^0)^2 (1 - (c_k^0 \tilde c_k^0)^2)
\sim \sum_{k=1}^{r} t_k^2 \alpha_k \beta_k (1 - (c_k^0 \tilde c_k^0)^2)
\end{align}

We now turn to the weighted estimator $\what \X_0 = \Omega \what \X \Pi^T$. From the weighted orthogonality condition, $\what \X = \what \X^{\dd}$, the optimal diagonal denoiser. From Theorem \ref{thm-diag}, the AMSE of $\what \X_0$ may be written
\begin{align}
\label{amse-X0}
\| \what \X_0 - \X_0\|_{\Fr}^2 
= \| \Omega ( \what \X^{\shr} - \X )\|_{\Fr}^2 
\sim \sum_{k=1}^{r} t_k^2 \alpha_k \beta_k \left( 1 
    - c_k^2 \tilde c_k^2 \cdot \frac{\alpha_k}{c_k^2\alpha_k + s_k^2 \mu} 
    \cdot \frac{\beta_k}{\tilde c_k^2\beta_k + \tilde s_k^2 \nu}\right).
\end{align}

Comparing \eqref{amse-Xshr} and \eqref{amse-X0}, the result will be proven if we can show
\begin{align}
\label{ineq6700}
(c_k^0)^2 < c_k^2  \cdot \frac{\alpha_k}{c_k^2\alpha_k + s_k^2 \mu}, 
\quad 1 \le k \le r,
\end{align}
and
\begin{align}
(\tilde c_k^0)^2 
< \tilde c_k^2  \cdot \frac{\beta_k}{\tilde c_k^2\beta_k + \tilde s_k^2 \nu}, 
\quad 1 \le k \le r.
\end{align}

By the symmetry in the problem, it is enough to prove \eqref{ineq6700}. Because we are working with each singular component separately, we will drop the subscript $k$.
From Proposition \ref{prop-vanilla}, the formula for $(c^0)^2$ is given by
\begin{align}
(c^0)^2 = 
\begin{cases}
\frac{1 - \gamma_0 / (t^0)^4}{1 + \gamma_0 / (t^0)^2},
        & \text{ if } t^0  > \gamma_0^{1/4}, \\
0, & \text{ if } t^0 \le \gamma_0^{1/4}.
\end{cases}
\end{align}

If $t^0 \le \gamma_0^{1/4}$, then \eqref{ineq6700} is trivial. Consequently, we assume $t^0 > \gamma_0^{1/4}$. Because $t^0 = t \sqrt{\alpha \beta / \nu}$ and $\gamma_0 = \gamma \mu / \nu$, this is equivalent to the condition
\begin{align}
t^4 > \gamma \frac{\mu \nu}{\alpha^2 \beta^2}.
\end{align}

Defining $R = \gamma / t^4$, we may consequently assume that
\begin{align}
R < \frac{\alpha^2 \beta^2}{\mu \nu} \le 1.
\end{align}

We may rewrite $(c^0)^2$ in terms of $R$ as follows:
\begin{align}
\label{eq6500}
(c^0)^2 = \frac{\alpha^2 \beta^2 - R \mu \nu}{\alpha^2 \beta^2 + t^2 \alpha \beta \mu R}.
\end{align}

From the formula
\begin{align}
c^2 = \frac{1 - \gamma / t^4}{1 + \gamma/t^2},
\end{align}
we may rewrite the right side of \eqref{ineq6700} as
\begin{align}
\label{eq7900}
\frac{c^2 \alpha}{c^2 \alpha + s^2 \mu}
= \frac{\alpha^2 \beta^2 - R \mu \nu + R(\mu\nu - \alpha^2 \beta^2)}
    {\alpha^2 \beta^2 + t^2 \alpha \beta \mu R 
        + \alpha \beta R[t^2 \mu \beta + \mu \beta - \alpha \beta - t^2\mu]}.
\end{align}

Comparing \eqref{eq6500} to \eqref{eq7900}, the inequality \eqref{ineq6700} is equivalent to showing
\begin{align}
\label{ineq6900}
(\alpha^2 \beta^2 - R\mu \nu)[t^2 \mu \beta + \mu \beta - \alpha \beta - t^2\mu]
< (\mu \nu - \alpha^2 \beta^2) (\alpha \beta + t^2 \mu R).
\end{align}

Because each side is affine linear in $R$, and $0 \le R \le \alpha^2 \beta^2 / (\mu \nu)$, it is enough to verify \eqref{ineq6900} at $R=0$ and $R = \alpha^2 \beta^2 / (\mu \nu)$. When $R = \alpha^2 \beta^2 / (\mu \nu)$, the left side of \eqref{ineq6900} is $0$, whereas the right side is non-negative because $\alpha < \sqrt{\mu}$ and $\beta < \sqrt{\nu}$, verifying the inequality in this case. When $R=0$, the difference between the right side and left side of \eqref{ineq6900}, divided by $\alpha \beta \mu$, is equal to
\begin{align}
\nu - \alpha^2 \beta^2 / \mu
    - \alpha \beta[t^2 \beta + \beta - \alpha \beta / \mu - t^2]
&= \nu - \alpha \beta[t^2 \beta + \beta - t^2]
\nonumber \\
&= \nu - \alpha \beta^2 t^2 - \alpha \beta^2 + \alpha \beta t^2
\nonumber \\
&=  t^2 \alpha \beta(1-\beta)+  \nu - \alpha \beta^2.
\end{align}
Since $\beta^2 \le \nu \le 1$ and $\alpha \le 1$, this expression is positive, verifying \eqref{ineq6900} and completing the proof.

\section{Proof of Proposition \ref{prop-estimating}}
\label{appendix-estimating}

From a standard Bernstein-type inequality for subexponential random variables (e.g.\ Proposition 5.16 from \cite{vershynin2010intro}), for every $\theta > 0$ we have:
\begin{align}
\Prob( \max_{i} |\hat a_i - \EE[\hat a_i] | > \theta) \le C p \exp\{-C'n \min(\theta^2 , \theta)\},
\end{align}
for constants $C, C' > 0$. Since $p \sim \gamma n$, the right hand side is summable in $n$; it follows from the Borel-Cantelli Lemma that $\max_{i} |\hat a_i - \EE[\hat a_i] | \to 0$ almost surely as $n \to \infty$. From the delocalization of $\X$'s singular vectors,  $\sup_{i,j} |X_{ij}|^2 = o(1/n)$. Since $\tr{\T}/n = 1$, we then have
\begin{align}
\EE[\hat a_i] = \EE\left[\sum_{j=1}^{n} X_{ij}^2 
    + 2\sum_{j=1}^{n} X_{ij} \sqrt{a_i b_j} G_{ij} +  \sum_{j=1}^{n} a_i b_jG_{ij}^2 \right]
= a_i + o(1).
\end{align}
Consequently, $\|\what \S_p - \S_p\|_{\op} = \max_{1 \le i \le p}|\what a_i - a_i| \to 0$ almost surely, as desired. Similar reasoning also shows that $\sum_{i=1}^{p} (\hat a_i - a_i)/n \to 0$ almost surely as $n \to \infty$. A nearly identical argument applied to the numerator of $\hat b_j$ then shows $\|\what \T_n - \T_n\|_{\op} \to 0$ almost surely, completing the proof.

\section{Proof of Proposition \ref{prop-snr}}
\label{appendix-snr}
To prove Proposition \ref{prop-snr},we begin by deriving a lower bound on the operator norm of the noise matrix $\N = \S^{1/2} \G \T^{1/2}$. We let $\a$ and $\b$ be unit vectors so that $ \G^T \b = \|\G\|_{\op} \a$. Then
\begin{align}
\|\G \T^{1/2}\|_{\op} \ge \| \T^{1/2} \G^T \b\| 
    = \|\G\|_{\op} \|\T^{1/2} \a\|.
\end{align}
%
%%%Since $\tilde{N}$ is an iid Gaussian matrix, the vector $a$ is uniformly distributed over the sphere in $\mathbb{R}^p$. Consequently, $\|B^{1/2} a\| \sim \sqrt{\tr{B} / p} \sim \sqrt{\E_{\nu \sim H_B}[ \nu]}$.

Next, we take unit vectors $\c$ and $\d$ so that $\G \T^{1/2}\d = \|\G \T^{1/2}\|_{\op} \c$. Then we have 
\begin{align}
\|\N\|_{\op}^2 \ge \|\S^{1/2} \G \T^{1/2} \d\|^2 
    = \|\G \T^{1/2}\|_{\op}^2 \|\S^{1/2}\c\|^2 
   \ge \|\G\|_{\op}^2 \cdot \|\T^{1/2} \a\|^2 \cdot \|\S^{1/2}\c\|^2 .
\end{align}
Since the distribution of $\G$ is orthogonally-invariant, the distributions of $\a$ and $\c$ are uniform over the unit spheres in $\R^n$ and $\R^p$, respectively. Consequently, $\|\T^{1/2} \a\|^2 \sim \tr{\T} / n$ and $\|\S^{1/2} \c\|^2 \sim \tr{\S} / p$. Therefore,
\begin{align}
\|\N\|_{\op}^2 \gtrsim (\tr{\S} / p) \cdot (\tr{\T} / n) \cdot \| \G \|_{\op}^2
\end{align}
where the inequality holds almost surely in the large $p$, large $n$ limit. Note that $\|\G\|_{\op} \sim 1 + \sqrt{\gamma}$ (see, e.g., \cite{bai2009spectral}), though we do not need to use this fact. 

Furthermore, we also have
\begin{align}
\tilde t_k^2 = t_k^2 \| \S^{-1/2} \u_k \|^2 \| \T^{-1/2} \v_k \|^2
\sim t_k^2 \cdot \frac{1}{p} \tr{\S^{-1}}\cdot \frac{1}{n} \tr{\T^{-1}}.
\end{align}
Consequently,
\begin{align}
\snr_k = \frac{t_k^2}{\|\N\|_{\op}^2} 
    &\lesssim \frac{t_k^2}{ (\tr{\S} / p) \cdot (\tr{\T} / n) \cdot  \| \G \|_{\op}^2 }
    \nonumber \\
&= \frac{t_k^2 \cdot (\tr{\S^{-1}} / p) \cdot (\tr{\T^{-1}} / n) }
                {\tau   \| \G \|_{\op}^2 }
= \frac{\tilde t_k^2}{\tau \| \G \|_{\op}^2}
    =  \frac{1}{\tau}\wtilde \snr_k,
\end{align}
completing the proof.

\section{Proof of Proposition \ref{prop-missing}}
\label{appendix-missing}

Let $\delta_{ij}$ be $1$ if entry $(i,j)$ is sampled, and $0$ otherwise. Then $\delta_{ij} \sim \text{Bernoulli}(p_i q_j)$. Let $\Delta = (\delta_{ij})$; then $\F^* (\F(\X)) = \Delta \odot \X$, where $\odot$ denotes the Hadamard product. Let $\q_r=(q_1^r,\dots,q_p^r)^T$ and $\q_c=(q_1^c,\dots,q_n^c)^T$.  The matrix $\Delta - \q_r \q_c^T$ is a random matrix with mean zero, whose entries are uniformly bounded. It follows from Corollary 2.3.5 of \cite{tao2012topics} that $\|\Delta - \q_r \q_c^T\|_{\op} / \sqrt{n} \le A$ a.s.\ as $n \to \infty$, for some constant $A > 0$.

We may write $\P \X \Q = \X \odot (\q_r \q_c^T)$, and consequently $\Delta \odot \X - \P \X \Q = (\Delta - \q_r \q_c^T) \odot \X$. Since $\X = \sum_{k=1}^{r} t_k \u_k \v_k^T$, it is enough to show that
\begin{align}
\| (\Delta - \q_r \q_c^T) \odot \u_k \v_k^T \|_{\op} \to 0
\end{align}
almost surely, for each $k$.

Suppose $\a$ and $\b$ are two unit vectors. Then 
\begin{align}
|\a^T [(\Delta - \q_r \q_c^T) \odot \u_k \v_k^T  ] \b|
= |(\a \odot \u_k)^T (\Delta - \q_r \q_c^T) (\b \odot \v_k)|
\le A \sqrt{n} \| \a \odot \u_k \| \| \b \odot \v_k\|
\end{align}
almost surely as $n \to \infty$.

Now, since $\a$ is a unit vector,
\begin{align}
\| \a \odot \u_k \| = \sqrt{\sum_{j=1}^{p} a_j^2 u_{jk}^2} \le \|\u_k\|_\infty
\end{align}
and similarly,
\begin{align}
\| \b \odot \v_k \|  \le \|\v_k\|_\infty.
\end{align}

Since $\max_{1\le k \le r}\|\u_k\|_\infty \|\v_k\|_\infty = o(n^{-1/2})$, the result follows.

\end{document}